\documentclass[VANCOUVER,STIX1COL]{WileyNJD-v2}

\articletype{Article Type}%

\received{26 April 2016}
\revised{6 June 2016}
\accepted{6 June 2016}

\usepackage{graphics}
\usepackage{amsmath}

\newcommand{\bbb}{\mbox{\boldmath $b$}}

\newcommand{\eee}{\mbox{\boldmath $e$}}
\newcommand{\fff}{\mbox{\boldmath $f$}}

\newcommand{\ppp}{\mbox{\boldmath $p$}}
\newcommand{\qqq}{\mbox{\boldmath $q$}}
\newcommand{\rrr}{\mbox{\boldmath $r$}}

\newcommand{\ttt}{\mbox{\boldmath $t$}}
\newcommand{\uuu}{\mbox{\boldmath $u$}}
\newcommand{\vvv}{\mbox{\boldmath $v$}}

\newcommand{\xxx}{\mbox{\boldmath $x$}}
\newcommand{\yyy}{\mbox{\boldmath $y$}}
\newcommand{\zzz}{\mbox{\boldmath $z$}}

\newcommand{\bmm}[1]{\mbox{\boldmath $ #1 $}}

\newcommand{\beq}{\begin{equation}}
\newcommand{\eeq}[1]{\label{#1} \end{equation}}
\newcommand{\beqa}{\begin{eqnarray}}
\newcommand{\eeqa}[1]{\label{#1} \end{eqnarray}}
\newcommand{\bmat}[1]{\left ( \begin{array}{#1}}
\newcommand{\emat}{\end{array} \right )}

\usepackage{amsmath}

\begin{document}

%\title{Convergence Acceleration of Preconditioned CG Solver Based on Error Vector Sampling for a Sequence of Linear Systems\protect\thanks{This is an example for title footnote.}}

\title{Convergence Acceleration of Preconditioned CG Solver Based on Error Vector Sampling for a Sequence of Linear Systems}

%\title{Effevtive Use of Error Vector Sampling: Convergence Acceleration of Preconditioned CG Solver  for a Sequence of Linear Systems and Condition Number Estimation}

\author[1]{Takeshi Iwashita*}

\author[2]{Kota Ikehara}

\author[1]{Takeshi Fukaya}

\author[3]{Takeshi Mifune}

\authormark{AUTHOR ONE \textsc{et al}}

\address[1]{\orgdiv{Information Initiative Center}, \orgname{Hokkaido University}, \orgaddress{Sapporo, \country{Japan}}}

\address[2]{\orgdiv{Graduate School of Information Science and Technology}, \orgname{Hokkaido University}, \orgaddress{Sapporo, \country{Japan}}}

\address[3]{\orgdiv{Graduate School of Engineering}, \orgname{Kyoto University}, \orgaddress{Kyoto, \country{Japan}}}

\corres{*Takeshi Iwashita, N 11 W 5, Sapporo, Japan. \email{iwashita@iic.hokudai.ac.jp}}

%\presentaddress{This is sample for present address text this is sample for present address text}

\abstract[Summary]{
In this paper, we focus on solving a sequence of linear systems with an identical (or similar) coefficient matrix. 
%The linear systems are sequentially processed due to the dependence of the right-hand side vector on the solution vector of the prior linear system. 
For this type of problems, we investigate the subspace correction and deflation methods, which use an auxiliary matrix (subspace) to accelerate the convergence of the iterative method. In practical simulations, these acceleration methods typically work well when the range of the auxiliary matrix contains eigenspaces corresponding to small eigenvalues of the coefficient matrix.
We have developed a new algebraic auxiliary matrix construction method based on error vector sampling, in which eigenvectors with small eigenvalues are efficiently identified in a solution process.
The generated auxiliary matrix is used for the convergence acceleration in the following solution step.
Numerical tests confirm that both subspace correction and deflation methods with the auxiliary matrix can accelerate the solution process of the iterative solver. 
Furthermore, we examine the applicability of our technique to the estimation of the condition number of the coefficient matrix.
The algorithm of preconditioned conjugate gradient (PCG) method with the condition number estimation is also shown.
}


%\jnlcitation{\cname{%
%\author{Williams K.}, 
%\author{B. Hoskins}, 
%\author{R. Lee}, 
%\author{G. Masato}, and 
%\author{T. Woollings}} (\cyear{2016}), 
%\ctitle{A regime analysis of Atlantic winter jet variability applied to evaluate HadGEM3-GC2}, \cjournal{Q.J.R. Meteorol. Soc.}, \cvol{2017;00:1--6}.}

\maketitle

\footnotetext{\textbf{Abbreviations:} CG, Conjugate Gradient; ICCG, Incomplete Cholesky Conjugate Gradient}

\section{Introduction}
A preconditioned conjugate gradient (CG) solver is widely used to solve a linear system of equations of a symmetric positive-definite (s.p.d.) matrix arising in various applications.
The computational time to solution is mostly given by the product of the number of iterations for convergence and the computational time per iteration.
Whereas high performance and parallel computing techniques are effective to reduce the computational time per iteration, the convergence acceleration of the solver is also an important topic.
It is well known that the convergence rate of the CG solver is affected by the condition number or the eigenvalue distribution of the coefficient matrix.
In practical simulations, the coefficient matrix often has a few small isolated eigenvalues, which lead to a significant decline in convergence. For these problems, the subspace correction~\cite{xu} and deflation~\cite{deflation0} methods are widely used to improve the convergence rate of the iterative solver.
%recognized as a convergence acceleration methods for this type of problems.

The procedure of the subspace correction and the deflation involves an auxiliary matrix to specify a certain subspace, in which errors are efficiently removed.
Therefore, a proper setting of the auxiliary matrix (subspace) is a key to make these acceleration methods work well.
For example, when the range of the matrix contains the eigenspaces corresponding to small isolated eigenvalues, the convergence rate of the solver is expected to be improved by the acceleration methods.
However, it is not easy to identify these eigenspaces.
Accordingly, in practical simulations, an effective auxiliary matrix is often derived from the knowledge of the problem.
For example, the coarse grid correction in the multigrid method~\cite{multigrid1, multigrid2}, which is regarded as one of the most successful subspace correction methods, uses the characteristics of discretized PDE problems.
Other examples of the auxiliary matrix or the subspace which is determined based on physics or models can be seen in the literature~\cite{deflation1, deflation2, deflation3, mifune1, igarashi-def} .
However, there are many cases in which the eigenvector with a small eigenvalue is hardly identified from the knowledge of the problem. 
For these problems, an automatic (algebraic) auxiliary matrix construction method that does not use special knowledge of the problem has been investigated. 

In this paper, we introduce an algebraic auxiliary matrix construction method for a problem involving a sequence of linear systems to be solved.
When the coefficient matrices are identical, it is often called a multiple right-hand side problem. 
In our method, an auxiliary matrix to specify the subspace is constructed using the sampling of error vectors in the preceding iterative solution process.
The idea is based on the expectation that the error which is not efficiently removed in the solution process contains useful information for the eigenvectors associated with small eigenvalues~\cite{JSIAM}.
% or  slowly convergent vectors~\cite{JSIAM}.
Although the error vector sampling during the solution process may seem difficult, it can be implemented by sampling the approximate solution vectors for the targeted problem.
After the solution process completes, the corresponding error vectors can be easily calculated.
We apply the Rayleigh-Ritz method using the subspace spanned by these error vectors to obtain (approximate) eigenvectors associated with small eigenvalues.
In our technique, the sampling plays a key role to save the additional memory footprint and computations for the subspace correction and the deflation, which is essential for many practical applications.

In this paragraph, we describe related works on algebraic auxiliary matrix construction for the convergence acceleration methods. 
Many related works can be found in the context of recycling Krylov subspace, deflation, augmented Krylov subspace, subspace recycling, and spectral preconditioning. 
After some early activities on the deflation in a GMRES solver~\cite{deflation-GMRES, deflation-GMRES2}, Moorgan proposed the GMRES-DR method.
In the method, basis vectors generated in the Arnoldi process in a restart period are used to determine the subspace for deflation~\cite{GMRES-DR}.
Moorgan et al. also introduced some variants of the GMRES-DR method which includes an application to the flexible GMRES method~\cite{GMRES-DR2, FGMRES-DEFLATION}.
Carpenter describes five major methods to specify the subspace (enrichment vectors) in the context of solvers based on the GMRES method~\cite{NASA}.
For CG solvers, Saad et al. introduced the deflated Lanczos algorithm and developed the deflated-CG method~\cite{saad-cg}. In this method, the vectors (subspace) used for the deflation are based on $A$-orthogonal basis vectors and are updated in the multiple linear system solution steps.
Abdel-Rehim et al. introduced the deflated restarted Lanczos algorithm~\cite{restart-lan}.
The techniques mentioned above were enhanced for nonlinear application problems, for example, in the research~\cite{Misha, Gosselet}.
As a recently published work, we refer to the paper~\cite{Daas}, in which Daas et. al. introduced a method based on the singular value decomposition. 
Moreover, it is noted that a block Krylov method can be used together with the convergence acceleration methods, though it is a popular technique for a multiple right-hand side problem in itself~\cite{block-GMRES}.
Finally, we refer to a recent survey paper written by Soodhalter et. al~\cite{survay-subspace}. The paper gives a comprehensive review of subspace recycling techniques to possibly cover most of related works to our research.

To the best of our knowledge, the above related papers do not explicitly discuss our approach based on error vector sampling.
In this paper, we describe the auxiliary matrix construction method based on vector sampling for the subspace preconditioning and the deflation method. We also introduce a cost model for the convergence acceleration.
Finally, we report the numerical results using test matrices of various application areas which were derived from the SuiteSparse Matrix Collection~\cite{Florida}, though our preliminary analyses only dealt with two computational electromagnetic problems~\cite{magn}.
The numerical results confirm the effectiveness of our method in terms of the convergence (\# iterations) and the computational time.
The numerical test also verifies our cost model and shows how the small eigenvalues are captured.
Furthermore, we show that our method can be used for the condition number estimation without significant additional computations in the iterative solution process.

\section{Problem definition}
In this paper, we deal with solving a sequence of $n$-dimensional linear systems:
\beq
A_{k} \xxx_{k}=\bbb_{k}, (k=1, 2, \ldots, k_{t}), 
\eeq{ak2}
where the coefficient matrix $A_{k}$ is a real symmetric positive-definite matrix.
We assume that the right-hand side vector $\bbb_k$ depends on the previous solution vectors. 
Consequently, the linear systems are solved sequentially. 
In this paper, we discuss the case where the coefficient matrices are all identical; 
\beq
A_{k} = A, (k=1, 2, \ldots, k_{t}). 
\eeq{ak3}
However, the technique introduced in the following sections is expected to work when the coefficient matrix changes but not dramatically. More precisely,
when the coefficient matrices have identical eigenvectors associated with small eigenvalues, it is possibly effective. 
In this paper, we solve the linear system of equations (\ref{ak2}) using a preconditioned Conjugate Gradient (CG) solver.

\section{Convergence Acceleration for Iterative Linear Solvers}
\subsection{Convergence Acceleration Methods}
In an iterative linear solver, its convergence rate directly affects the solution time. 
%One of the most popular method is preconditioning, which is essential for practical problems. 
%In the method, a linear system is transformed into a preconditioned system which is better conditioned to attain a better convergence rate.
In this paper, we focus on convergence acceleration methods that use a (user-specified) subspace different from the subspace designated by the coefficient matrix such as the Krylov subspace.
In these methods, the dimension of the subspace used is typically much smaller than $n$, and the error component involved in the subspace is efficiently removed by a special procedure.
A multigrid method can be regarded as a typical example of this type of convergence acceleration method.
In this paper, we discuss the subspace correction and deflation methods, both of which use a user-specified subspace to accelerate the convergence.

\subsection{Subspace Correction Method}
The subspace correction is a generalized version of the coarse grid correction of the multigrid method. We describe its procedure for an $n$-dimensional linear system; $A\xxx=\bbb$, where $\xxx$ is the unknown vector, and $\bbb$ is the right-hand-side vector. 

In the subspace correction method, an approximate solution vector $\tilde{\xxx}$ is updated as follows:
\begin{description}
\item{Step 1:} Compute $\fff=W^{\top} (\bbb-A\tilde{\xxx})$
\item{Step 2:} Solve $ (W^{\top} AW)\uuu=\fff$
\item{Step 3:} Update $\tilde{\xxx} \gets \tilde{\xxx}+W\uuu$
\end{description}
$W$ is the auxiliary matrix to designate the user-specified subspace.
The number of columns of $W$ is typically much less than $n$.
%In a practical situation, the approximation vector is often updated by the smoothing process before Step 1 and after Step 3. The smoothing process is typically given by one or a few stationary iterations.

When we use the subspace correction method together with a Krylov subspace method, we construct the preconditioner based on the correction like the multigrid (2-level) preconditioning~\cite{multigrid1}.
The subspace correction preconditioning\footnote{In this paper, we use the word ``subspace correction preconditioning'', which appears in the references~\cite{SC1, SC2}. The preconditioning based on the same concept is often called 2-level preconditioning, or spectral preconditioning~\cite{2level}, especially when the subspace are associated with eigenspaces.} can be combined with any other (standard) preconditioning techniques in the additive/multiplicative Schwarz preconditioning manner. When the stand-alone preconditioner is denoted by $M^{-1}$, the additive Schwarz subspace correction preconditioner  $M_{sc}^{-1}$ is given by
\beq
M_{sc}^{-1}=M^{-1}+W(W^{\top}AW)^{-1} W^{\top}.
\eeq{mc}
When the subspace preconditioning is only used, $M$ is given by the identity matrix $I$.

\subsection{Deflation method}
In this subsection, we describe the procedure of the deflated CG method~\cite{saad-cg} for $A \xxx = \bbb$.
In the deflation method, we use the projector given by
\begin{equation}
\label{deflation:projector}
P = I - W (W^\top A W)^{-1} (A W)^\top.
\end{equation}
$P$ decomposes the $n$-dimensional space $\mathbb{R}^{n}$ into two $A$-orthogonal spaces $\mathcal{W}$ and $\mathcal{W}^{\bot}$.
By using the projector, the solution vector $\xxx$ can be split into two components:
\beq
\xxx = \yyy + \zzz, \ \yyy=(I-P) \xxx, \ \zzz=P \xxx.
\eeq{split}
In the deflation method, two vector components $\yyy$ and $\zzz$ are individually derived.
The vector $\yyy$ is in the lower dimensional space $range(W)$ and is given by
\beq
\yyy = (I-P) \xxx = W (W^\top A W)^{-1} W^\top \bbb.
\eeq{x1}
Because it holds that $P^{\top} A (I -P) = O$, the second component $\zzz$ is computed by solving the deflated system
\beq
P^{\top} A \zzz = P^{\top} \bbb.
\eeq{x2}
In this paper, the deflated system having a semi-positive definite coefficient matrix (\ref{x2}) is solved using a preconditioned CG solver.
Algorithm \ref{alg:DPCG} shows the algorithm of the deflated CG method. It is noted that the projector $P$ is not explicitly constructed in practical implementations.

\begin{algorithm}[t]
\caption{Deflated PCG method}
\label{alg:DPCG} 
\begin{algorithmic}[1]
\Require $A$, $\bmm{b}$, $M$, $W$, $P$, $\bmm{x}_{0}$, $\varepsilon$
\State $\bmm{r}_0 = P^\top ( \bmm{b} - A \bmm{x}_0 )$
\State $\bmm{p}_0 = M^{-1} \bmm{r}_0$
\For{$i = 0, 1, 2, \ldots $ \textbf{until} $ \| \bmm{r}_i \|_2 \le \varepsilon \| \bmm{b} \|_2$}
\State $\displaystyle \alpha_i = \frac{(M^{-1}\bmm{r}_{i}, \bmm{r}_{i})}{(\bmm{p}_i, P^\top A \bmm{p}_i)}$
\State $\bmm{x}_{i+1} = \bmm{x}_i + \alpha_i \bmm{p}_i$
\State $\bmm{r}_{i+1} = \bmm{r}_i - \alpha_i P^\top A \bmm{p}_i$
\State $\displaystyle \beta_i = - \frac{(M^{-1}\bmm{r}_{i+1}, \bmm{r}_{i+1})}{(M^{-1}\bmm{r}_{i}, \bmm{r}_{i})}$
\State $\bmm{p}_{i+1} = M^{-1} \bmm{r}_{i+1} + \beta_i \bmm{p}_i$
\EndFor
\State $\bmm{x} = P \bmm{x}_i + W (W^\top A W)^{-1} W^\top \bmm{b}$
\Ensure $\bmm{x}$
\end{algorithmic}
\end{algorithm}

%\begin{algorithm}[t]
%\caption{Deflated PCG method}
%\label{alg:DPCG} 
%\begin{algorithmic}[1]
%\Require $\bmm{A}$, $\bmm{b}$, $\bmm{M}$, $\bmm{B}$, $\bmm{P}$, $\bmm{x}_{0}$, $\varepsilon$
%\State $\bmm{r}_0 = \bmm{P}^\top ( \bmm{b} - \bmm{A} \bmm{x}_0 )$
%\State $\bmm{p}_0 = \bmm{M}^{-1} \bmm{r}_0$
%\For{$i = 0, 1, 2, \ldots $ \textbf{until} $ \| \bmm{r}_i \|_2 \le \varepsilon \| \bmm{b} \|_2$}
%\State $\displaystyle \alpha_i = \frac{(\bmm{M}^{-1}\bmm{r}_{i}, \bmm{r}_{i})}{(\bmm{p}_i, \bmm{P}^\top \bmm{A} \bmm{p}_i)}$
%\State $\bmm{x}_{i+1} = \bmm{x}_i + \alpha_i \bmm{p}_i$
%\State $\bmm{r}_{i+1} = \bmm{r}_i - \alpha_i \bmm{P}^\top \bmm{A} \bmm{p}_i$
%\State $\displaystyle \beta_i = - \frac{(\bmm{M}^{-1}\bmm{r}_{i+1}, \bmm{r}_{i+1})}{(\bmm{M}^{-1}\bmm{r}_{i}, \bmm{r}_{i})}$
%\State $\bmm{p}_{i+1} =  \bmm{M}^{-1} \bmm{r}_{i+1} + \beta_i \bmm{p}_i$
%\EndFor
%\State $\bmm{x} = \bmm{P} \bmm{x}_i + \bmm{W} (\bmm{W}^\top \bmm{A} \bmm{W})^{-1} \bmm{W}^\top \bmm{b}$
%%
%%
%\Ensure $\bmm{x}$
%\end{algorithmic}
%\end{algorithm}

\section{Auxiliary matrix construction method based on error vector sampling}
\subsection{Auxiliary matrix based on eigenvectors}
In the subspace correction and deflation methods, the key to the convergence acceleration is in proper setting of the auxiliary matrix $W$. 
Typically, when the range of $W$ contains eigenspaces corresponding to small eigenvalues of the coefficient matrix, the methods work.
In practical problems, a coefficient matrix often has a few isolated small eigenvalues, which worsens the convergence of iterative solver. 
These eigenvalues typically arises from the physical property of the targeted problem.

Let us consider the situation that $W$ is an $n \times 1$ matrix and $W = [ \uuu_{s} ]$, where $\uuu_{s}$ is the eigenvector associated with the smallest eigenvalue $\lambda_{s}$.
%
%$W$ consists of the eigenvector $\uuu_{s}$ associated with the smallest eigenvalue $\lambda_{s}$. 
%The auxiliary matrix $W$ is an $n \times 1$ matrix.
%
We assume that $\lambda_{s} \ll 1$ and it is isolated.
It is also assumed that the coefficient matrix has an eigenvalue close to or larger than 1.
In this case, the subspace correction preconditioning $M_{sc}^{-1}$ with $M=I$ only shifts the eigenvalue $\lambda_{s}$ to $\lambda_{s}$+1.
In other words, the preconditioned coefficient matrix has an eigenvalue of $\lambda_{s}$+1 and $n-1$ eigenvalues that are identical to those of $A$ and larger than $\lambda_{s}$.
Consequently, the condition number of the preconditioned coefficient matrix is better than that of $A$, which results in better convergence for the preconditioned system. 

When we use the deflation method with the above setting for $W$,  $\lambda_{s}$ is removed in the coefficient matrix of (\ref{x2}), $P^{\top} A$. 
$P^{\top} A$ has a zero eigenvalue which is associated with $\uuu_{s}$, and other eigenvalues and eigenvectors are the same as $A$.
The (preconditioned) CG method can be applied to ($\ref{x2}$) because $P^{\top} \bbb$ is involved in range($P^{\top} A$), and its convergence rate is improved from that for the original linear system, $A \xxx = \bbb$.

The above discussion is straightforwardly extended to the case that $W$ consists of multiple eigenvectors associated with small eigenvalues.
However, calculation of eigenvalues and eigenvectors typically requires more computational efforts than solving the linear system itself.
Consequently, in practical simulations, the knowledge of the problem is often used for identifying the eigenvectors associated with small eigenvalues and constructing a proper auxiliary matrix. 
But, there are problems in which the origin of the small eigenvalue is unclear from the viewpoint of physics or simulation models.
In this paper, we focus on a problem of solving a sequence of linear systems, and intend to develop an automatic auxiliary matrix construction method for the problem.

\subsection{Auxiliary matrix construction method based on error vector sampling}
This subsection describes our auxiliary matrix construction method based on error vector sampling for a sequence of linear systems (\ref{ak2}).
During the first iterative solution process for $A \xxx_{1}=\bbb_{1}$, we preserve $m$ approximation solution vectors $\tilde{\xxx}_1^{(s)}  (s=1, 2, \ldots, m )$.
Typically, $m$ is much smaller than $n$.
After the solution process is completed, the error vectors that correspond to $\tilde{\xxx}_1^{(s)}$ are calculated by 
\beq
\eee^{(s)}=\xxx_1 - \tilde{\xxx}_1^{(s)} \, (s=1, 2, \ldots, m ).
\eeq{ek}
Applying the Gram--Schmidt process to these error vectors, we obtain the mutually orthogonal $\bar{m}$ $(\leq$$\, m)$ normal basis vectors:
\beq
\bar{\eee}^{(1)}, \, \bar{\eee}^{(2)}, \, \cdots, \, \bar{\eee}^{(\bar{m})}.  
\eeq{ekd}
In our technique, the Rayleigh--Ritz method based on the space spanned by $\bar{\eee}^{(s)}$ is used to identify approximate eigenvectors associated with small eigenvalues of $A$.

The auxiliary matrix construction method is given as follows:

Step 1: Solve the $\bar{m}$-dimensional eigenvalue problem~\footnote{In this paper, we intend to identify eigenvectors with relatively small eigenvalues of the coefficient matrix itself. But, it is possible to consider identifying eigenvectors with relatively small eigenvalues of the preconditioned matrix. In this case, we should use the preconditioned matrix instead of $A$ in (\ref{Rayleigh-Ritz}).}:
\beq
E^{\top} A E \ttt = \lambda \ttt, 
\eeq{Rayleigh-Ritz}

where
\beq
E = [\bar{\eee}^{(1)} \ \bar{\eee}^{(2)} \ \cdots \ \bar{\eee}^{(\bar{m})}] 
\eeq{Rayleigh-Ritz2}

Step 2: When Ritz value $\lambda$ is less than a preset threshold $\theta$,  Ritz vector $E \ttt$ is adopted as a column vector of $W$. The number of Ritz values less than $\theta$ is denoted by $\tilde{m}$, and the Ritz vector that corresponds to each small Ritz value is written as $E\ttt_i  (i=1, 2, \ldots, \tilde{m})$. Finally, the auxiliary matrix $W$ is given by
\beq
W=[E\ttt_1 \ E\ttt_2 \ \cdots \ E\ttt_{\tilde{m}}].
\eeq{bk}
The threshold is typically much less than 1, namely $(0<\theta \ll 1)$ when the coefficient matrix is diagonally (or properly) scaled.

\subsection{Selection Method for Stored Approximation Vectors}
In practical analyses, to avoid an excessive additional cost (in memory space and computations), the number of stored vectors, $m$, should be substantially small. We use the selection method based on ``sampling''. We intend to store approximate solution vectors with a certain interval in the solution process. 
Considering the difficulty of prediction of the number of iterations for convergence, we use following two methods for sampling.
In the sampling method A, we use the algorithm shown in Appendix A.
When $m$ is set to be 4 and the (preconditioned) CG solver attains convergence at 1,000-th iteration, the sampling method preserve the approximation vectors at 256, 384, 512, and 768-th iterations.
%
%The interval of the sampling in the method is given by $2^{d}$ or $2^{d+1}$, where $d$ is a positive integer that satisfies $d<\mbox{log}_2 (N/m)\leq d+1$, and $N$ is the total iteration count.
%
Another method (sampling method B) is based on the relative residual norm.
We take a sample of approximation vector when the relative residual norm first reaches $10^{-s\alpha/(m+1)},  (s=1, 2, \ldots, m)$, when the convergence criteria is given by $10^{-\alpha}$.
Based on the preliminary test result, we use the sampling method A when we do not explicitly mention the sampling method.

% and $l_{max}$ is a given parameter that satisfies $m^{l_{max}}>N_{max}$, where $N_{max}$ is the maximum iteration count of the linear solver. If the iterative solver converges at the 5,000-th iteration (i.e., $N$=5,000) with the setting of $m$=4, the approximate solution vectors at the 1,024, 2,048, 3,072, and 4,096-th iterations are preserved. 

%Although it is difficult to predict the number of iterations necessary for convergence, Algorithm 1 can preserve $m$ approximation solution vectors that satisfy the following property with minimum memory space, where $\tilde{\xxx}_{i_{te}}$ is the approximate solution vector at the $i_{te}$-th iteration of the linear solver. The interval in the iteration count between the vectors is $2^{d}$ or $2^{d+1}$, where $d$ is a positive integer that satisfies $d<\mbox{log}_2 (N/m)\leq d+1$. In the algorithm, $N$ is the total iteration count and $l_{max}$ is a given parameter that satisfies $m^{l_{max}}>N_{max}$, where $N_{max}$ is the maximum iteration count of the linear solver. If the iterative solver converges at the 5,000-th iteration (i.e., $N$=5,000) with the setting of $m$=4, the approximate solution vectors at the 1,024, 2,048, 3,072, and 4,096-th iterations are preserved. 

\subsection{Computational Cost for Subspace Correction Preconditioning and Deflation}
In this subsection, we discuss the additional computational cost for two convergence acceleration techniques. 
%\subsubsection{Subspace correction preconditioning}
%
Computational time per iteration of preconditioned CG solver $T$ is split into two parts:
\beq
T=T_{pre}+T_{cg},
\eeq{T} 
where $T_{pre}$ and $T_{cg}$ are the computational time for preconditioning and CG solver parts, respectively.
Because the total data amount for matrices and vectors are typically larger than cache memory in practical simulations, most of computational kernels of the solver results in memory-bound. Consequently, we estimate the computational time using the amount of transferred data from main memory.
In the analysis, double precision floating point numbers are used for matrices and vectors.
The main part of the CG solver is a sparse matrix vector multiplication (spMV) kernel. 
We estimate the amount of transferred data for spMV as $20n+12nnz$, where $nnz$ is the number of nonzero elements of $A$ and the unit is Byte.
Although the cache hit ratio for elements of the source vector depends on the nonzero pattern of $A$, we use relatively optimistic estimation.
The transferred data for other parts that consist of inner products and vector updates is estimated as $56n$.
When the effective memory bandwidth is denoted by $b_{m}$ Byte/s, $T_{cg}$ is estimated as
\beq
T_{cg}=(76n+12nnz)/b_{m}.
\eeq{tcg}
When we use IC preconditioning, the transferred data for preconditioning is almost the same as spMV. 
Finally, the computational time for an ICCG iteration that is denoted by $T_{iccg}$ is approximately given by
\beq
T_{iccg}=(100n+24nnz)/b_{m}.
\eeq{ticcg}

When we consider the subspace correction (SC) preconditioning, the additional cost for 
$W(W^{\top}AW)^{-1} W^{\top}$ should be taken into account.
In the estimation, we ignore the cost for $(W^{\top}AW)^{-1}$ because the dimension $\tilde{m}$ is much smaller than $n$ on the setting of $m \ll n$.
%In other words, the number of sampling vectors, $m$, must be much smaller than $n$.
The additional transferred data for the SC preconditioning is mainly for the $n \times \tilde{m}$ dense matrix $W$, and it is estimated as 
$16\tilde{m}n+16n$.
When we use SC preconditioning together with IC preconditioning, the computational time for a SC-ICCG iteration that is denoted by $T_{sciccg}$ 
is estimated as
\beq
T_{sciccg}=(116n+16\tilde{m}n+24nnz)/b_{m}.
\eeq{tsciccg}
From (\ref{ticcg}) and (\ref{tsciccg}), we can (roughly) estimate the ratio of the computational cost per iteration for two solvers, SC-ICCG and ICCG, which is denoted by $\gamma_{sciccg}$, as follows:
\beq
\gamma_{sciccg}=(116+16\tilde{m}+24nnz_{av})/(100+24nnz_{av}),
\eeq{ganma}
where $nnz_{av}$ is the average number of nonzero elements per row.
When the number of iteration of SC-ICCG is less than $1/\gamma_{sciccg}$ of that of ICCG, SC-ICCG is expected to outperform ICCG. 

Next, we consider the deflation method.
When we use the deflation method, the additional cost is in calculating $P^{\top}A$.
The data transferred for $P^{\top}A$ is estimated to be almost the same as SC preconditioning because both $AW$ and $W$ are dense matrices with the identical size. Consequently, (\ref{ganma}) can be used for the ICCG solver with deflation.

Based on the expectation in the reduction of the iteration count and (\ref{ganma}), we can set the number of sample vectors, $m$.
For example, when we expect a 40\% reduction by the convergence acceleration method for a problem of $nnz_{av}$$=\,$30, $\tilde{m}$ $(\leq$$\,m)$ should be less than $20$.

%
%\begin{algorithm}                      
%\caption{Selection of approximate solution vectors}         
%\label{alg1}                          
%\begin{algorithmic}                  
%\State $h=1$
%\For{$i_{te}=1, 2, \cdots$}
%\State {\bf Solver part}
%\State Convergence check
%\If  {$(\mbox{mod}(i_{te}, h)==0)$} 
%\State  $i_{t}=\sum_{l=0}^{l_{max}} (-1)^{l} \lfloor (i_{te}-1)/m^{l} \rfloor $
%\State  $j=\mbox{mod}(i_{t}, m)+1$
%\State  $\tilde{\xxx}^{(j)}=\tilde{\xxx}_{i_{te}}$
%\If {$(i_{te}==h*m)$} 
%\State $h=h*2$ 
%\EndIf
%\EndIf
%\EndFor
%\end{algorithmic}
%\end{algorithm}

\section{Numerical Results}
\subsection{Test Conditions}
Numerical tests were conducted to examine the effect of convergence acceleration methods (subspace correction and deflation) based on our algebraic auxiliary matrix generation method.
For the test matrix, we downloaded 30 relatively large matrices from the SuiteSparse Matrix Collection~\cite{Florida} and applied the diagonal scaling to them.
We picked up symmetric positive definite matrices that were mainly derived from computational science or engineering problems.
Table \ref{matInfo} shows properties of the test matrices.
For each coefficient matrix, we solve a linear system of equations 6 times.
The convergence criterion is given by the relative residual 2-norm being less than $10^{-8}$. 
After the first solution process is completed, the auxiliary matrix is generated and used in the following 5 solution processes, in which the solver performance is evaluated.
For the right-hand vector, we used two kinds of vectors; a vector of ones and a random vector.
In the former case, we solve an identical linear system 6 times. 
When we use random vectors, each linear system to be solved is different to one another.
In this paper, we report the result when the number of sampled vectors, $m$, is set to be 20.

Numerical tests were conducted on a computational node of Fujitsu CX2550 (M4) at Information Initiative Center, Hokkaido University.
The node is equipped with two Intel Xeon (Gold6148, Skylake) processors, each of which has 20 cores, and 384GB memory.
The program code was written in C and OpenMP for the thread parallelization. 
Intel C compiler version 19.1.3.304 was used with the option of ``-O3 -qopenmp -ipo -xCORE-AVX512''.
In the tests of parallel multithreaded solvers, 40 threads were used.

\begin{table*}[tbp]
%      \small
	\centering
%	\vspace{0.75\baselineskip}
	\caption{Matrix information for the test problems}
	\label{matInfo}
%	\begin{tabular}{|c|c|c|c|c|}
	\begin{tabular}{llrrr}
	\hline
	Data set & Problem type & Dimension & \# nonzero & $nnz_{av}$\\
	\hline
Queen\_4147& 2D/3D Problem& 4,147,110& 316,548,962& 76.3 \\
Bump\_2911& 2D/3D Problem& 2,911,419& 127,729,899& 43.9 \\
G3\_circuit& Circuit Simulation Problem& 1,585,478& 7,660,826& 4.8 \\
Flan\_1565& Structural Problem& 1,564,794& 114,165,372& 73.0 \\
Hook\_1498& Structural Problem& 1,498,023& 59,374,451& 40.0 \\
StocF-1465& Computational Fluid Dynamics Problem& 1,465,137& 21,005,389& 14.3 \\
Geo\_1438& Structural Problem& 1,437,960& 60,236,322& 41.9 \\
Serena& Structural Problem& 1,391,349& 64,131,971& 46.1 \\
thermal2& Thermal Problem& 1,228,045& 8,580,313& 7.0 \\
ecology2& 2D/3D Problem& 999,999& 4,995,991& 5.0 \\
bone010& Model Reduction Problem& 986,703& 47,851,783& 48.5 \\
ldoor& Structural Problem& 952,203& 42,493,817& 44.6 \\
audikw\_1& Structural Problem& 943,695& 77,651,847& 82.3 \\
Emilia\_923& Structural Problem& 923,136& 40,373,538& 43.7\\
boneS10& Model Reduction Problem& 914,898& 40,878,708& 44.7 \\
PFlow\_742& 2D/3D Problem& 742,793& 37,138,461& 50.0 \\
tmt\_sym& Electromagnetics Problem& 726,713& 5,080,961& 7.0 \\
apache2& Structural Problem& 715,176& 4,817,870& 6.7 \\
Fault\_639& Structural Problem& 638,802& 27,245,944& 42.7\\
parabolic\_fem& Computational Fluid Dynamics Problem& 525,825& 3,674,625& 7.0 \\
bundle\_adj& Computer Vision Problem& 513,351& 20,207,907& 39.4 \\
af\_shell8& Subsequent Structural Problem& 504,855& 17,579,155& 34.8 \\
af\_shell4& Subsequent Structural Problem& 504,855& 17,562,051& 34.8 \\
af\_shell3& Subsequent Structural Problem& 504,855& 17,562,051& 34.8 \\
af\_shell7& Subsequent Structural Problem& 504,855& 17,579,155& 34.8 \\
inline\_1& Structural Problem& 503,712& 36,816,170& 73.1 \\
af\_0\_k101& Structural Problem& 503,625& 17,550,675& 34.8 \\
af\_4\_k101& Structural Problem& 503,625& 17,550,675& 34.8 \\
af\_3\_k101& Structural Problem& 503,625& 17,550,675& 34.8 \\
af\_2\_k101& Structural Problem& 503,625& 17,550,675& 34.8 \\ \hline
\end{tabular}
\end{table*}

\subsection{Numerical results on the sequential solver}
%\subsubsection{Tests using a right-hand side vector of ones}
\subsubsection{Performance evaluation}
\begin{figure*}[tbp]
\centering
\includegraphics[clip, scale=0.9,  bb=50 90 560 400]{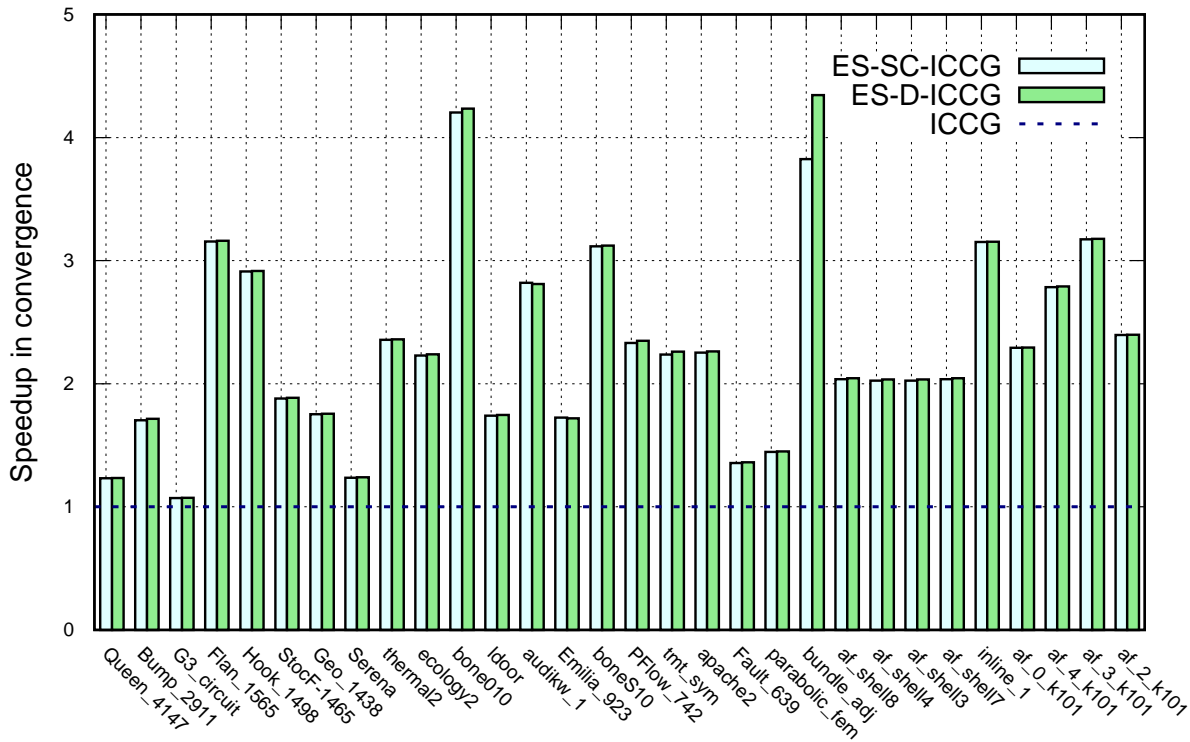}
\caption{Speedup in convergence of ES-SC-ICCG and ES-D-ICCG over ICCG (\bbb: random vector)}
\label{speedup-con-seq}
\end{figure*}

Table \ref{result-1-ones} lists the numerical results of the standard ICCG solver and its variants with the introduced convergence acceleration techniques, when a vector of ones is used for the right-hand side. 
ES-SC-ICCG denotes the CG solver with IC and subspace correction preconditioning based on the proposed error vector sampling method. ES-D-ICCG denotes the deflated ICCG solver using our technique. 
The table show the average computational time (sec) in 5 solution steps, which is denoted by $T_{t}$. 
Table \ref{result-1-rand} shows the results when random vectors are used for the right-hand side. The table shows the average number of iterations and computational time in 5 solution steps.
The numerical results indicate that both solvers based on the proposed method achieve convergence acceleration for all 60 test cases (30 datasets $\times$ 2 kinds of right-hand side vectors). 
The convergence acceleration was significant for some of datasets.
In the numerical tests using the vector of ones, the acceleration method attains more than 3-fold speedup in convergence for 16 out of 30 datasets.
Even when we used random right-hand side vectors, the convergence was more than twice as fast as the ICCG solver for 20 out of 30 datasets as shown in Fig. \ref{speedup-con-seq}.

\begin{figure}[tbp]
\centering
\includegraphics[clip, scale=0.7,  bb=50 50 500 350]{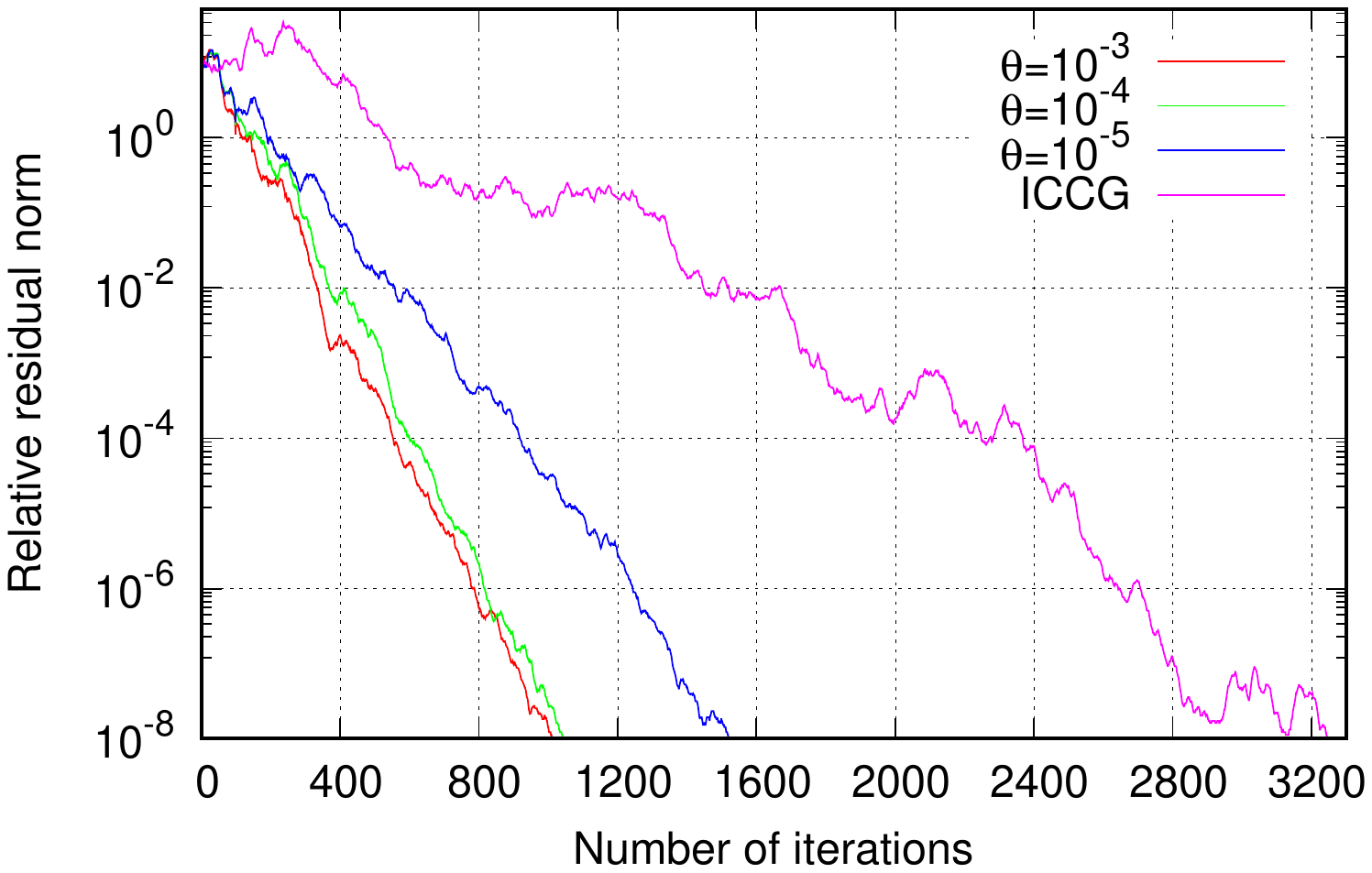}
\caption{Convergence behavior of ES-SC-ICCG (dataset: Flan\_1565)}
\label{flan-r-sc}
\end{figure}

\begin{figure}[tbp]
\centering
\includegraphics[clip, scale=0.7,  bb=50 50 500 350]{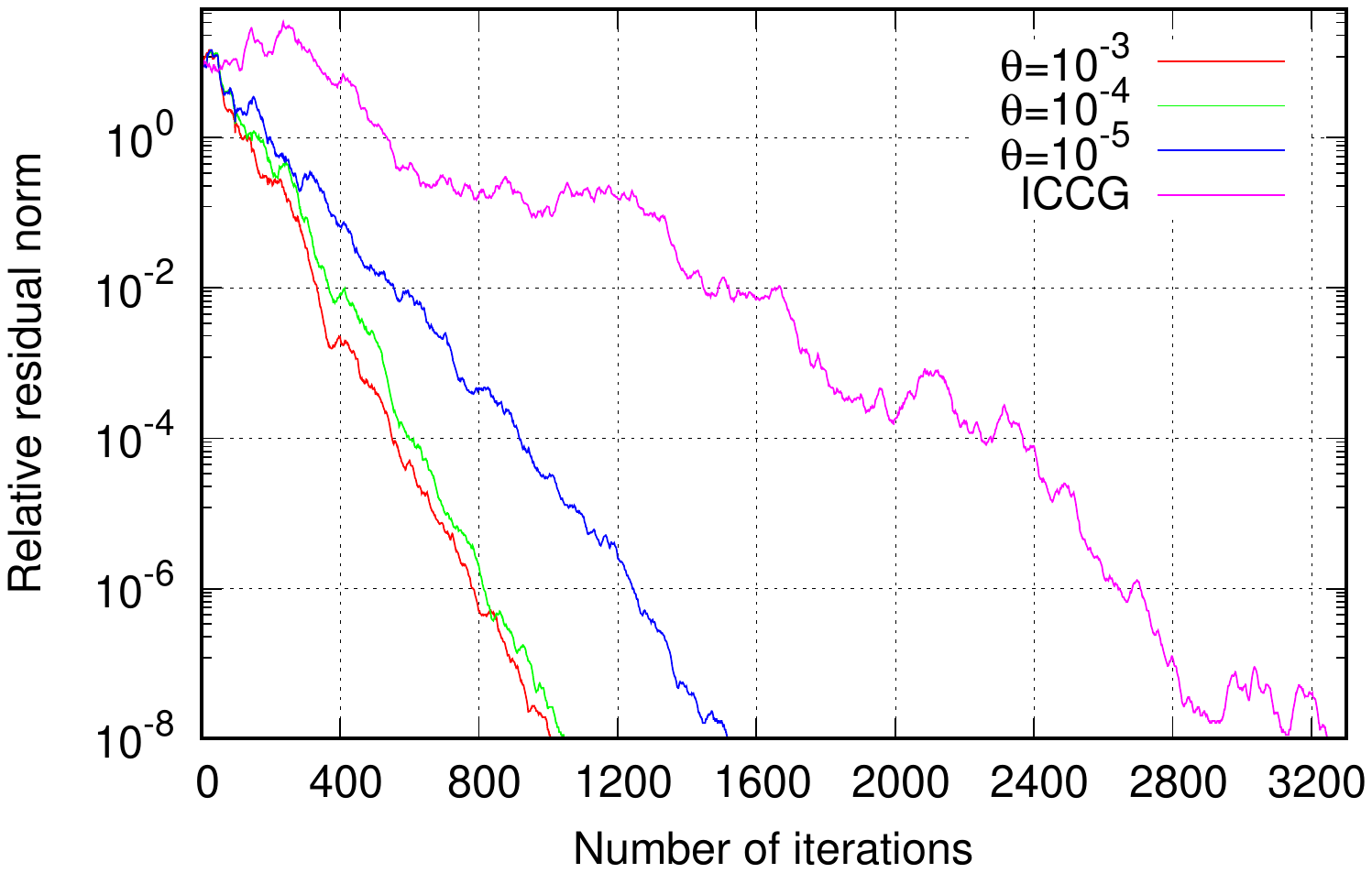}
\caption{Convergence behavior of ES-D-ICCG (dataset: Flan\_1565)}
\label{flan-r-d}
\end{figure}

\begin{figure}[tbp]
\centering
\includegraphics[clip, scale=0.7,  bb=50 50 500 350]{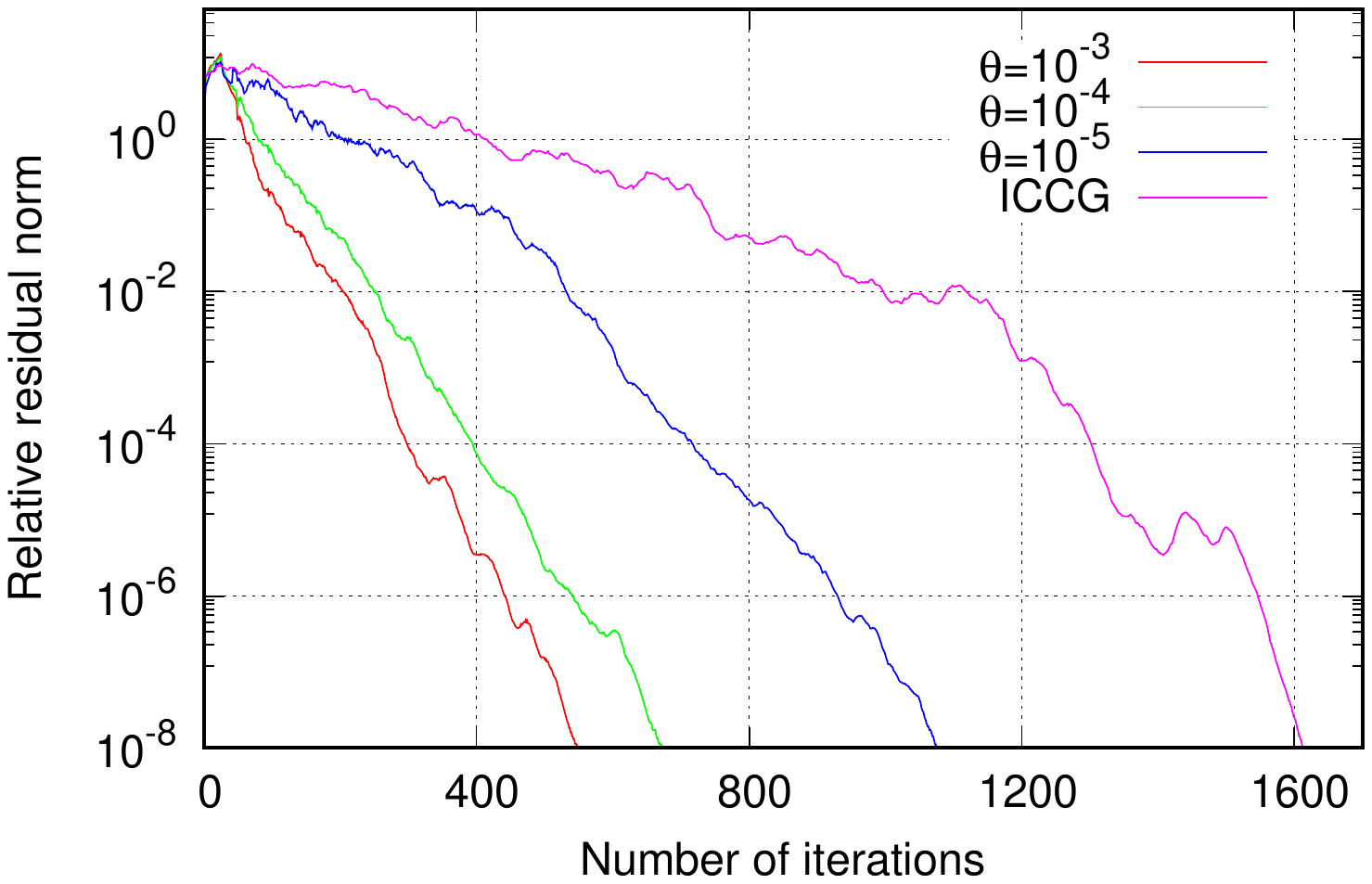}
\caption{Convergence behavior of ES-SC-ICCG (dataset: Hook\_1498)}
\label{hook-r-sc}
\end{figure}

\begin{figure}[tbp]
\centering
\includegraphics[clip, scale=0.7,  bb=50 50 500 350]{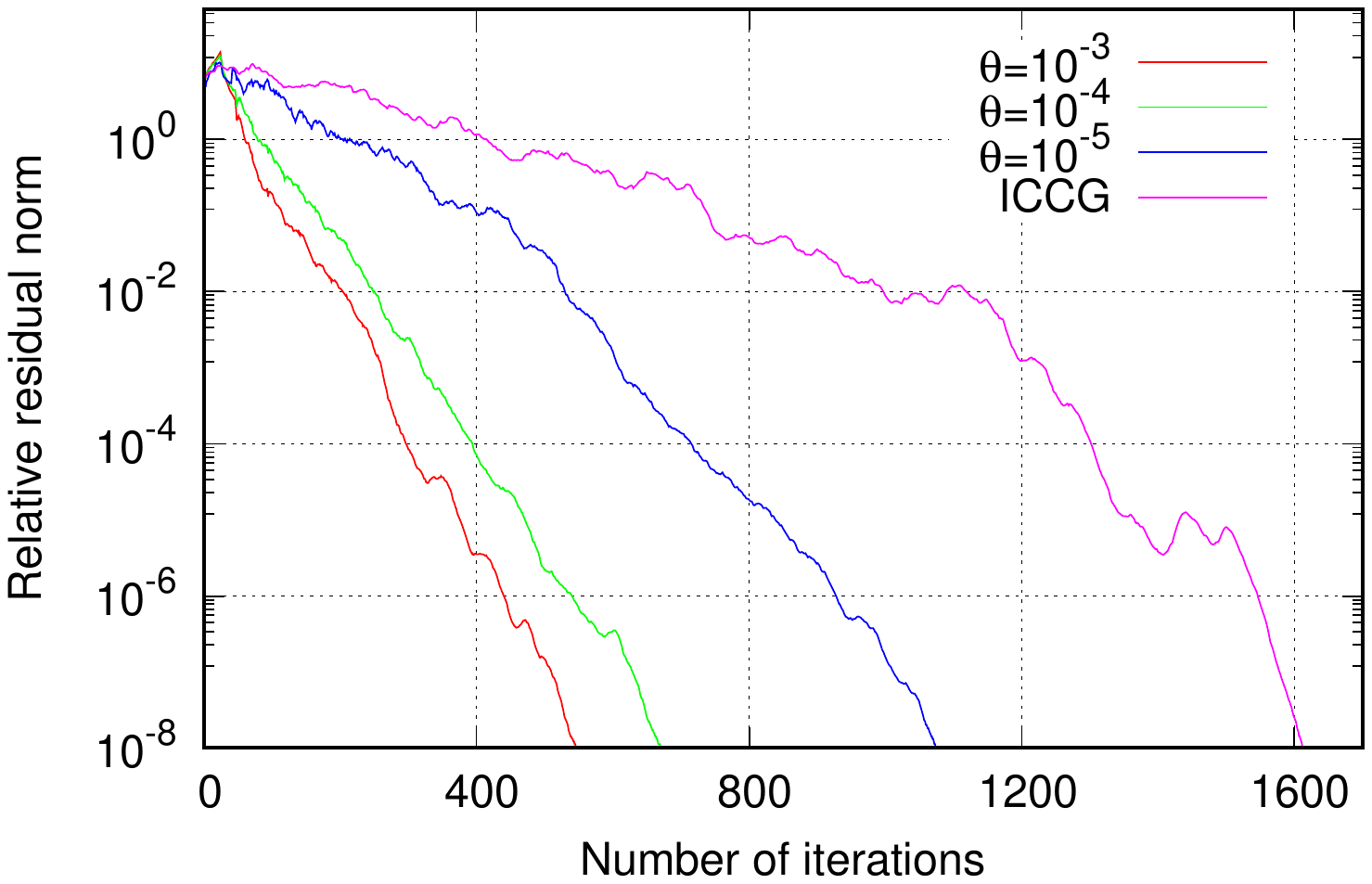}
\caption{Convergence behavior of ES-D-ICCG (dataset: Hook\_1498)}
\label{hook-r-d}
\end{figure}

Figures \ref{flan-r-sc}, \ref{flan-r-d}, \ref{hook-r-sc}, and \ref{hook-r-d} show the convergence behaviors of ES-SC-ICCG and ES-D-ICCG solvers for Flan\_1565 and Hook\_1498 datasets when a random vector is used for the right-hand side. The figures also confirm the effectiveness of the subspace correction and the deflation based on our technique. Numerical results imply that the larger $\theta$, which typically leads to larger $\tilde{m}$, results in the better convergence. This characteristic is also confirmed by the result listed in Tables \ref{result-1-ones} and \ref{result-1-rand}.  
Figures \ref{flan-r-sc}, \ref{flan-r-d}, \ref{hook-r-sc}, and \ref{hook-r-d} demonstrate that the convergence behaviors of two solvers are identical, though the treatments for the slow convergent error (eigenvectors associated with small eigenvalues) are different between two solvers~\cite{Zhao}.
While we examined the convergence behavior of the residual norm for all test cases, we observed that the convergence properties of two solvers were almost the same for most of test cases.
The result indicates that the effects of the SC preconditioning and the deflation are similar when the coefficient matrix is diagonally scaled and the identical subspace that corresponds to eigenvectors associated with small eigenvalues is used.
%Consequently, we may say that there is no practically significant difference in using either of the subspace correction and the deflation if the auxiliary matrix is identical, though ES-D-ICCG is more or less better than ES-SC-ICCG in most of test cases.

\begin{figure*}[tbp]
\centering
\includegraphics[clip, scale=0.9,  bb=50 90 560 330]{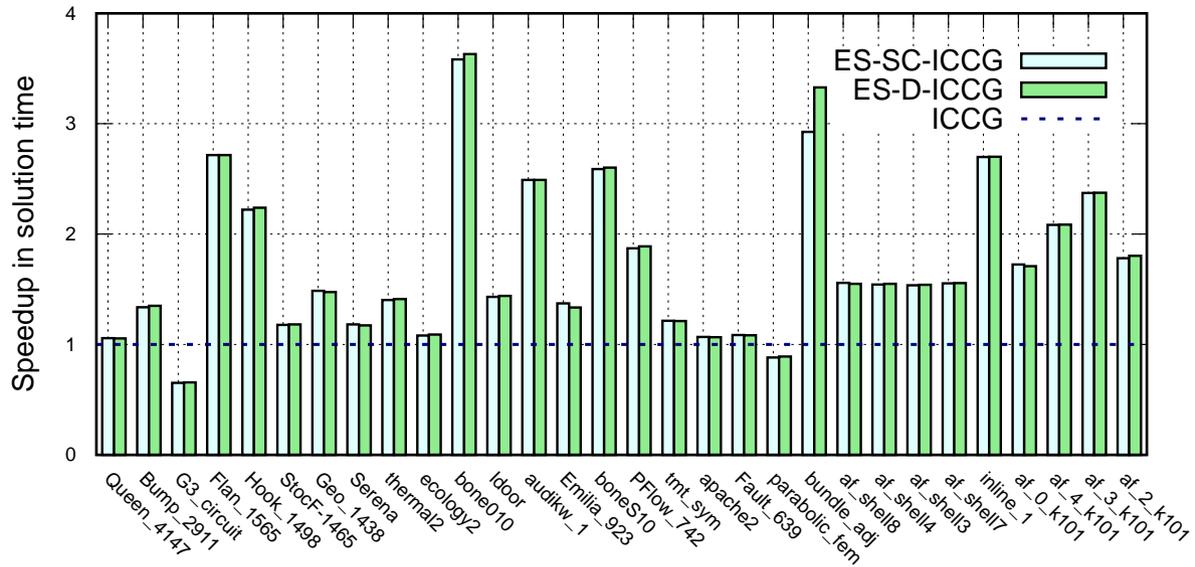}
\caption{Speedup in computational time of ES-SC-ICCG and ES-D-ICCG over ICCG (\bbb: random vector)}
\label{speedup-seq}
\end{figure*}

Next, we examine the computational time to solution. Table \ref{result-1-ones} shows that the solution time is reduced in 28 out of 30 cases in the tests using the right-hand side vector of ones.
For 16 datasets, the computational time of the solvers using our technique (ES-SC- and ES-D-ICCG) is reduced to less than half of that of the normal ICCG solver. The performance difference between two solvers, ES-SC-ICCG and ES-D-ICCG is marginal. 
In the numerical test using random vectors, the computational time is also reduced in 28 out of 30 cases.
Table \ref{result-1-rand}  and  Fig. \ref{speedup-seq} show the effectiveness of our technique in the random vector test.
In these tests, performance improvement is not attained in the G3\_circuit and parabolic\_fem datasets, which have relatively small $nnz_{av}$ values. In (\ref{ganma}), $\gamma_{sciccg}$ is enlarged when $nnz_{av}$ decreases.
It means that it becomes difficult to obtain performance improvement in the solution time by the subspace correction preconditioning and the deflation method. In other words, for a dataset with a small $nnz_{av}$ value, the convergence rate should be substantially improved by the limited number of sample vectors to achieve solver performance improvement.
In the numerical test, ES-SC-ICCG and ES-D-ICCG solvers obtained their best results for 12 out of 30 datasets when $\tilde{m}$ is equal to $m$ ($=$20).
For these datasets, an increase in the number of sample vectors, $m$, possibly improves the solver performance.
%
%Finally, Fig. \ref{speedup-seq} shows the speedup ratio of ES-SC-ICCG and ES-D-ICCG over ICCG, which confirms the effectiveness of our technique.

\subsubsection{Verification of the model for computational time per iteration}
\begin{figure}[tbp]
\centering
\includegraphics[clip, scale=0.7,  bb=50 50 500 350]{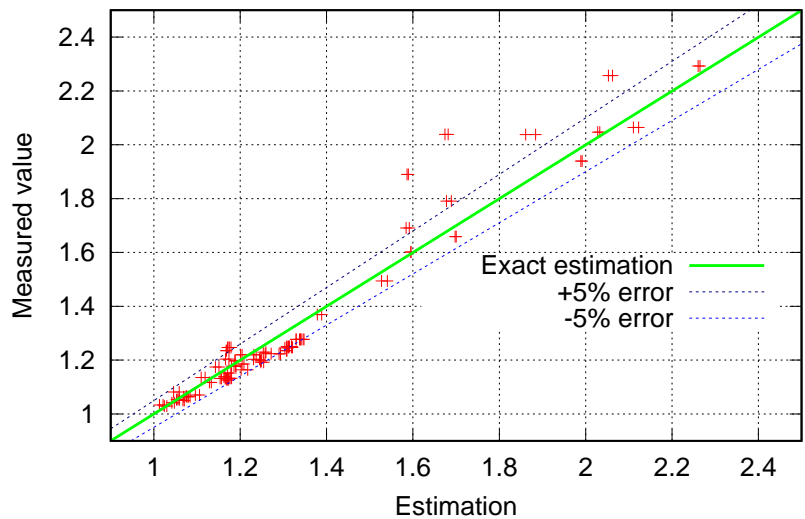}
\caption{Comparison of estimated and measured values of computational time per iteration}
\label{model}
\end{figure}

The application of subspace correction or deflation typically leads to increase in the computational cost per iteration. 
In this subsection, we examine the performance model for iteration cost introduced in Sec. 4.4.
In Fig. \ref{model}, we plot the measured and estimated values for the ratio of computational time of an ICCG iteration to that of an ES-SC- or ES-D-ICCG iteration. The estimated value for two solvers are given by (\ref{ganma}).
Figure \ref{model} shows the result for all test cases, though only one mark is plotted for an identical $\tilde{m}$.
For most of test cases, equation (\ref{ganma}) gives good estimation for the ratio, and the error of the estimation is within $\pm$5\%. Consequently, (\ref{ganma}) can be used for estimation of the additional cost for subspace correction or deflation.
However, in some test cases, especially when the measured value is over 2.0, relatively large estimation error is observed. These results arise in the G3\_circuit, ecology2, and apache2 datasets.
The coefficient matrices of these datasets commonly have a small number of nonzero elements per row ($nnz_{av}$) and a relatively structured nonzero element pattern. Namely, these matrices are derived from relatively simple problems, and (\ref{ticcg}) tends to give an overestimation for such a problem.  
Moreover, (\ref{ganma}) implies that the impact of the additional cost for the convergence acceleration on the computational time tends to be large when $nnz_{av}$ is small.
Accordingly, we recommend that the number of sampling vectors $m$ (the upper bound of $\tilde{m}$) should be small for a problem with small $nnz_{av}$.

\subsubsection{Other factors on solver performance}
{\it Sampling method}

In preliminary analyses, we compared two sampling methods A and B. 
Table \ref{sampling B} shows the results of the solver using the sampling method B for Flan\_1565 and Hook\_1498.
In comparison of tables \ref{result-1-ones} and  \ref{sampling B}, the sampling method A gives better convergence acceleration than the method B.
Because this tendency was observed for other test datasets, we decided to mainly use the sampling method A in our numerical tests.
Moreover, the numerical test implies that the additional sampling of the approximation vector when the residual norm increases or stagnates is effective for improvement of the convergence acceleration effect. Because it is not straightforward to mathematically interpret the phenomenon, we intend to investigate the behavior of the error in the solution process in our future work based on numerical tests.  
 
{\it Sampling of residual vectors}

In this paper, we consider the sampling of a relatively small number of vectors because it is practically important to save the additional memory space and computational cost. Considering other related techniques, the sampling of residual vectors might be of interest.
We have an intuitive perspective for the comparison of sampling of error vectors and residual vectors.
Because it holds that $A \eee_{s} = \rrr_{s}$,  the component corresponding small eigenvalues in $\eee_{s}$ is numerically reduced in $\rrr_{s}$ by multiplication of  $A$, where $\eee_{s}$ and $\rrr_{s}$ are the sampled error and residual vectors, respectively.
Consequently, it is expected that the error vector sampling is superior to the residual vector sampling to capture (approximate) eigenvectors that corresponds to small eigenvalues, which leads to better preconditioning effect for convergence.
To verify our perspective, we conducted additional numerical tests of the solver using the residual vector sampling.
In the numerical test using Flan\_1565 and Hook\_1498, it is shown that a small Ritz value less than $10^{-1}$ cannot be obtained and the convergence acceleration of the subspace correction and the deflation does not work well.
Considering the numerical results, we can say that the error vector sampling outperforms the residual vector sampling to construct an effective mapping operator for subspaces used in the convergence acceleration techniques.

{\it Verification of Ritz vector}

\begin{figure}[tbp]
\centering
\includegraphics[clip, scale=0.6,  bb=30 30 530 790]{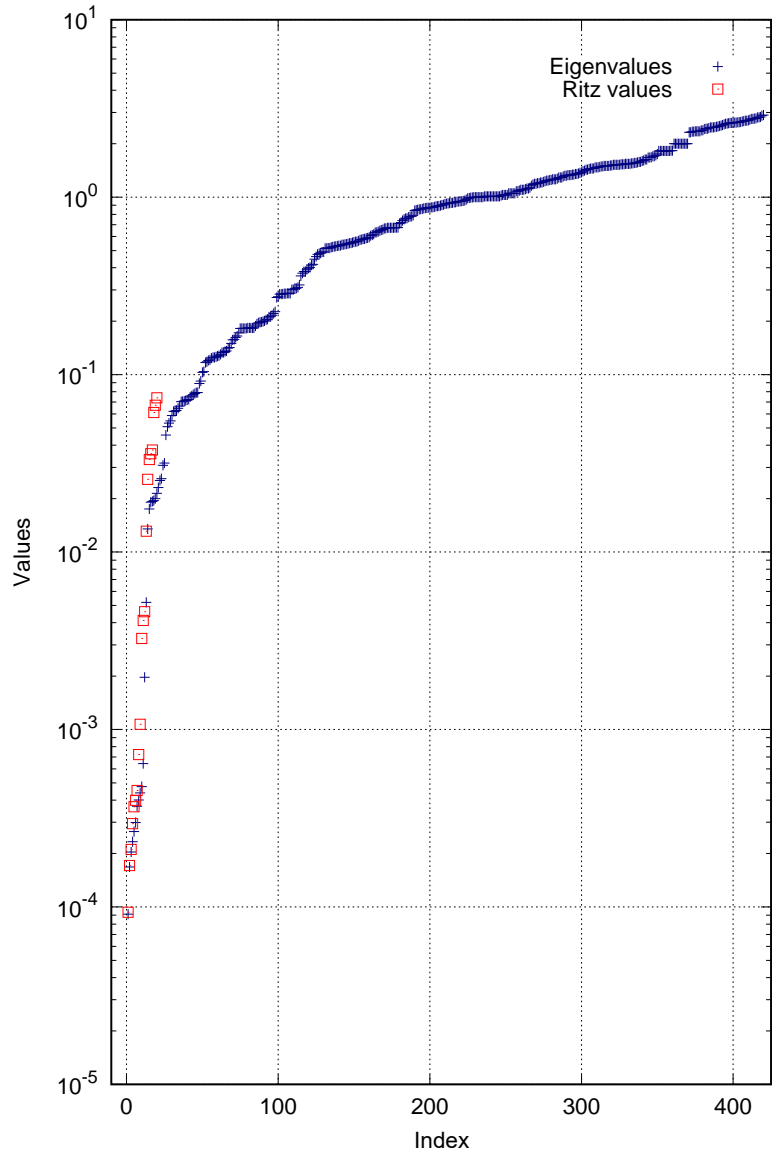}
\caption{Comparison of eigenvalues and Ritz values (dataset: bcsstik06, $n$=420, $m$=20)}
\label{eigen}
\end{figure}

\begin{figure}[tbp]
\centering
\includegraphics[clip, scale=0.7,  bb=50 50 500 350]{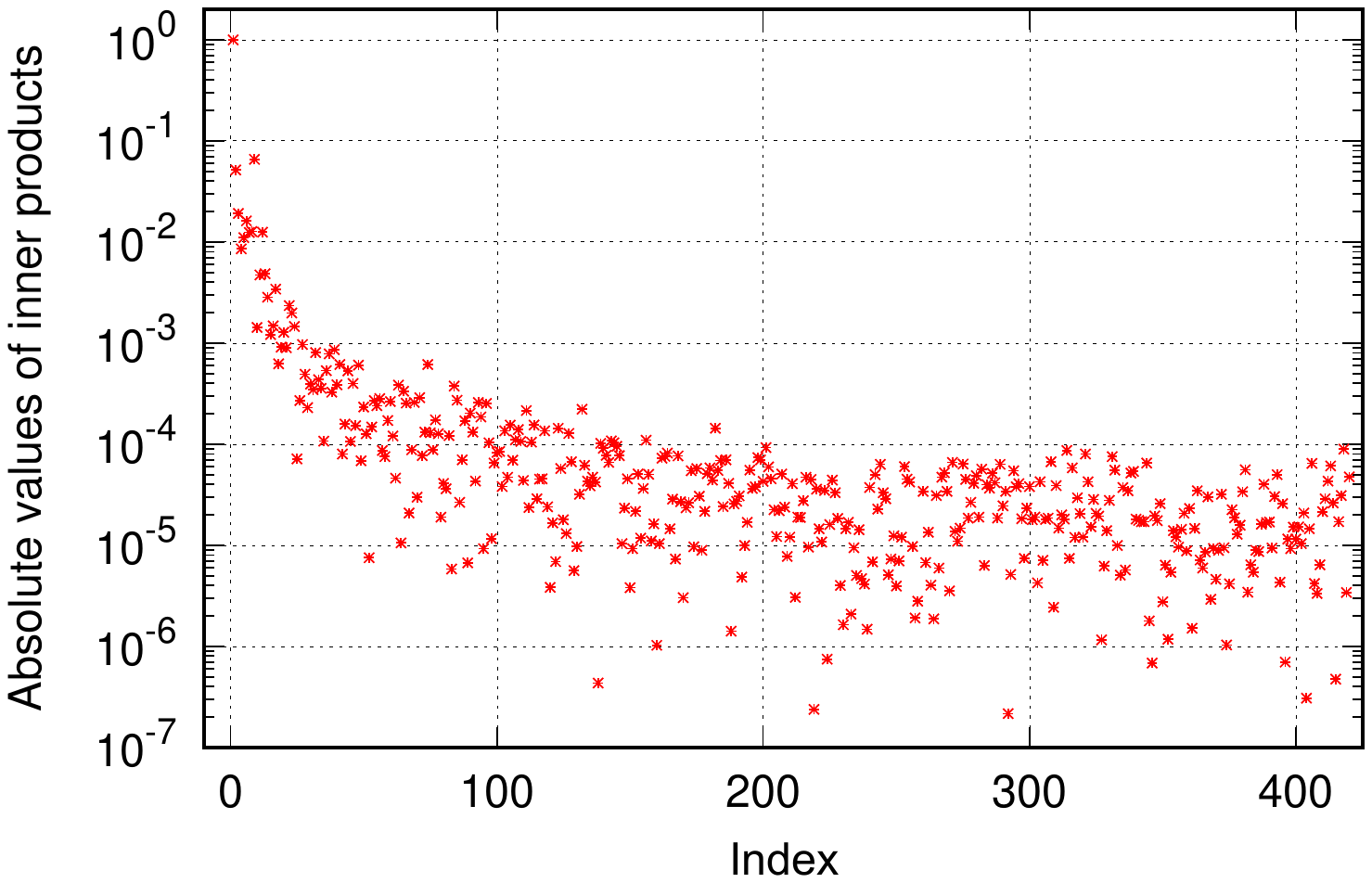}
\caption{The absolute values of inner products, $\| (\tilde{\vvv}_{1}, \vvv_{ir})\|, (ir=1, \ldots, 420)$}
\label{eigen2}
\end{figure}

In this section, we try to examine the property of Ritz vector calculated in our technique using a small size dataset (bccstk06: a 420 $\times$ 420 matrix).
Figure \ref{eigen} plots the eigenvalue distribution of the coefficient matrix and the Ritz values obtained in our method applied to a non-preconditioned CG solver.
It is confirmed that some small eigenvalues including the smallest eigenvalue are well approximated by the obtained Ritz values.
Moreover, we checked the orthogonality of the normalized Ritz vector that corresponds to the smallest Ritz value, $\tilde{\vvv}_{1}$, to the normalized eigenvectors of $A$ denoted by $\vvv_{ir}, (ir=1,\ldots,420)$. It is noted that $ir$ is the index of eigenvalues in ascending order. Figure \ref{eigen2} shows the absolute value of the inner product $(\vvv_{r}, \vvv_{ir})$.
The magnitude of $|(\tilde{\vvv}_{1}, \vvv_{1})|$ is close to 1 and it is substantially larger than those for other inner products, most of which are less than $10^{-3}$.

%\begin{table*}
%\centering
%(a) ICCG 
%
%\begin{tabular}{llllllllll} \hline
%\multicolumn{2}{l}{Queen\_4147} & \multicolumn{2}{l}{Bump\_2911} & \multicolumn{2}{l}{G3\_circuit} & \multicolumn{2}{l}{Flan\_1565} & \multicolumn{2}{l}{Hook\_1498} \\	\hline		
%\#Ite. & $T_{t}$ & \#Ite. & $T_{t}$ & \#Ite. & $T_{t}$ & \#Ite. & $T_{t}$ & \#Ite. & $T_{t}$ \\ \hline
%\end{tabular}
%
%%& StocF-1465 
%%Geo_1438			Serena			thermal2			ecology2		
%
%
%\end{table*}

\begin{table*}
%\normalsize
\fontsize{9pt}{0.4cm}\selectfont

%\small
\centering
\caption{Numerical results (sequential solver, $\bbb=(1, 1, \ldots, 1)^\top$)}
\label{result-1-ones}
\begin{tabular}{lllllllllllllllll} \hline
 & & \multicolumn{3}{l}{Queen\_4147} & \multicolumn{3}{l}{Bump\_2911} & \multicolumn{3}{l}{G3\_circuit} & \multicolumn{3}{l}{Flan\_1565} & \multicolumn{3}{l}{Hook\_1498} \\	\cline{3-17}
Solver & $\theta$ & $\tilde{m}$ & \#Ite. & $T_{t}$ & $\tilde{m}$ & \#Ite. & $T_{t}$ & $\tilde{m}$ & \#Ite. & $T_{t}$ & $\tilde{m}$ & \#Ite. & $T_{t}$ & $\tilde{m}$ & \#Ite. & $T_{t}$ \\ \hline
ICCG&&-&3128 &2763 &-&1551 &584 &-&898 &44.8 &-&3124 &996 &-&1617 &287 \\ \hline 
ES-SC-ICCG&$10^{-3}$&20&995 &1039 &20&526 &249 &18&705 &70.1 &20&1082 &398 &20&472 &108 \\ 
&$10^{-4}$&19&2041 &2121 &18&824 &382 &9&707 &54.8 &19&1212 &449 &13&676 &144 \\
&$10^{-5}$&7&2816 &2667 &5&1118 &445 &1&887 &49.6 &8&1766 &596 &5&1080 &209 \\ \hline
ES-D-ICCG&$10^{-3}$&20&993&1036 &20&459&218 &18&702&70.1 &20&942 &347 &20&469 &108 \\ 
&$10^{-4}$&19&2044&2120 &18&821&381 &9&706&54.1 &19&1213 &443 &13&675 &144 \\
&$10^{-5}$&7&2818 &2670 &5&1117 &449 &1&887 &49.8 &8&1762 &595 &5&1078 &209 \\ \hline
\end{tabular}
\vspace{1\baselineskip}

\begin{tabular}{lllllllllllllllll} \hline
 & & \multicolumn{3}{l}{StocF-1465} & \multicolumn{3}{l}{Geo\_1438} & \multicolumn{3}{l}{Serena} & \multicolumn{3}{l}{thermal2}& \multicolumn{3}{l}{ecology2} \\ \cline{3-17}
Solver &$\theta$&$\tilde{m}$&\#Ite.&$T_{t}$&$\tilde{m}$&\#Ite.&$T_{t}$&$\tilde{m}$&\#Ite.&$T_{t}$&$\tilde{m}$&\#Ite.&$T_{t}$&$\tilde{m}$&\#Ite.&$T_{t}$ \\ \hline
ICCG&&-&56109&4741 &-&443&79.6 &-&301&55.7 &-&2281&141 &-&1823&49.7  \\ \hline
ES-SC-ICCG&$10^{-3}$&20&14780&2011 &15&248&58.2 &7&243&49.6 &20&849&89 &20&813&50.1 \\
&$10^{-4}$&20&14775&2001 &2&387&72.5 &0&-&-&17&994&99 &15&902&48.5 \\
&$10^{-5}$&20&14775&1998 &0&-&-&0&-&-&4&1523&111 &5&1329&50.3 \\ \hline
ES-D-ICCG&$10^{-3}$&20&14731&1992 &15&248&55.9 &7&242&49.5 &20&847&89 &20&808&49.9 \\
&$10^{-4}$&20&14717&1988 &2&386&72.9 &0&-&-&17&992&99 &15&899&48.3 \\
&$10^{-5}$&20&14717&2001 &0&-&-&0&-&-&4&1519&111 &5&1328&49.7 \\ \hline
\end{tabular}
\vspace{1\baselineskip}

\begin{tabular}{lllllllllllllllll} \hline
 & & \multicolumn{3}{l}{bone010} & \multicolumn{3}{l}{ldoor} & \multicolumn{3}{l}{audikw\_1} & \multicolumn{3}{l}{Emilia\_923}  &\multicolumn{3}{l}{boneS10} \\ \cline{3-17}
Solver &$\theta$&$\tilde{m}$&\#Ite.&$T_{t}$&$\tilde{m}$&\#Ite.&$T_{t}$&$\tilde{m}$&\#Ite.&$T_{t}$&$\tilde{m}$&\#Ite.&$T_{t}$&$\tilde{m}$&\#Ite.&$T_{t}$ \\ \hline
ICCG&&-&4162 &801 &-&2160 &293 &-&2629 &583 &-&462 &53.6 &-&8532 &1275 \\ \hline
ES-SC-ICCG&$10^{-3}$&20&943 &213 &20&658 &111 &20&745 &185 &20&218 &32.2 &20&2688 &486 \\ 
&$10^{-4}$&18&967 &216 &16&1073 &174 &9&1138 &265 &19&266 &38.9 &20&2688 &487 \\
&$10^{-5}$&13&1302 &280 &3&1663 &238 &4&1521 &343 &5&373 &46.5 &20&2688 &488 \\ \hline
ES-D-ICCG&$10^{-3}$&20&935&211 &20&655&110 &20&756&188 &20&201&29.7 &20&2682 &485  \\
&$10^{-4}$&18&962&215 &16&1072&173 &9&1084&253 &19&265&38.9 &20&2682 &485 \\ 
&$10^{-5}$&13&1293 &278 &3&1662 &237 &4&1586 &358 &5&374 &46.8 &20&2682 &486 \\ \hline
\end{tabular}
\vspace{1\baselineskip}

\begin{tabular}{lllllllllllllllll} \hline
 & & \multicolumn{3}{l}{PFlow\_742} & \multicolumn{3}{l}{tmt\_sym}  & \multicolumn{3}{l}{apache2}  & \multicolumn{3}{l}{Fault\_639} & \multicolumn{3}{l}{parabolic\_fem} \\ \cline{3-17}
Solver &$\theta$&$\tilde{m}$&\#Ite.&$T_{t}$&$\tilde{m}$&\#Ite.&$T_{t}$&$\tilde{m}$&\#Ite.&$T_{t}$&$\tilde{m}$&\#Ite.&$T_{t}$&$\tilde{m}$&\#Ite.&$T_{t}$ \\ \hline
ICCG&&-&33076&3357 &-&1252&35.9 &-&768&16.5 &-&2187&177 &-&1131&18.9 \\ \hline
ES-SC-ICCG&$10^{-3}$&20&10359&1299 &20&507&27.0 &19&359&16.0 &20&806&83 &18&671&22.4 \\
&$10^{-4}$&20&10360&1299 &15&613&29.3 &12&429&15.7 &15&1366&134 &7&862&20.7 \\
&$10^{-5}$&20&10359&1299 &3&1013&34.7 &2&663&17.1 &4&1905&164 &0&-&- \\ \hline
ES-D-ICCG&$10^{-3}$&20&10269&1287 &20&501&26.5 &19&360&16.3 &20&798&82 &18&670&22.3 \\ 
&$10^{-4}$&20&10269&1287 &15&610&29.1 &12&428&15.8 &15&1364&133 &7&861&20.7 \\
&$10^{-5}$&20&10269&1285 &3&1011&33.9 &2&662&16.8 &4&1901&164 &0&-&- \\ \hline
\end{tabular}
\vspace{1\baselineskip}

\begin{tabular}{lllllllllllllllll} \hline
 & & \multicolumn{3}{l}{bundle\_adj} & \multicolumn{3}{l}{af\_shell8}  & \multicolumn{3}{l}{af\_shell4} & \multicolumn{3}{l}{af\_shell3} & \multicolumn{3}{l}{af\_shell7} \\ \cline{3-17}
Solver &$\theta$&$\tilde{m}$&\#Ite.&$T_{t}$&$\tilde{m}$&\#Ite.&$T_{t}$&$\tilde{m}$&\#Ite.&$T_{t}$&$\tilde{m}$&\#Ite.&$T_{t}$&$\tilde{m}$&\#Ite.&$T_{t}$ \\ \hline
ICCG&&-&42809 &2275 &-&1048 &52.0 &-&1048 &52.0 &-&1048 &52.3 &-&1048 &53.0 \\  \hline
ES-SC-ICCG&$10^{-3}$&20&11705 &824 &18&483 &31.4 &18&481 &31.1 &18&481 &31.4 &18&483 &31.5 \\ 
&$10^{-4}$&18&11533 &793 &9&614 &35.8 &9&615 &35.3 &9&615 &35.7 &9&614 &35.5 \\
&$10^{-5}$&17&11460 &781 &0&-&-&0&-&-&0&-&-&0&-&- \\ \hline
ES-D-ICCG&$10^{-3}$&20&9740 &686 &18&481 &31.4 &18&479 &31.0 &18&479 &31.4 &18&481 &31.5 \\ 
&$10^{-4}$&18&10117 &698 &9&613 &35.4 &9&615 &35.5 &9&615 &35.7 &9&613 &35.6 \\
&$10^{-5}$&17&10532 &717 &0&-&-&0&-&-&0&-&-&0&-&- \\ \hline
\end{tabular}
\vspace{1\baselineskip}

\begin{tabular}{lllllllllllllllll} \hline
 & & \multicolumn{3}{l}{inline\_1} & \multicolumn{3}{l}{af\_0\_k101}  & \multicolumn{3}{l}{af\_4\_k101}  & \multicolumn{3}{l}{af\_3\_k101}  & \multicolumn{3}{l}{af\_2\_k101} \\ \cline{3-17}
Solver &$\theta$&$\tilde{m}$&\#Ite.&$T_{t}$&$\tilde{m}$&\#Ite.&$T_{t}$&$\tilde{m}$&\#Ite.&$T_{t}$&$\tilde{m}$&\#Ite.&$T_{t}$&$\tilde{m}$&\#Ite.&$T_{t}$ \\ \hline
ICCG&&-&8487&879 &-&12953&636 &-&9993&489 &-&8519&423 &-&13092&648 \\ \hline
ES-SC-ICCG&$10^{-3}$&20&2573&311 &20&4153&276 &20&3093&204 &20&2632&176 &20&4194&279 \\
&$10^{-4}$&20&2572&310 &20&4153&276 &20&3094&204 &20&2633&176 &20&4194&279 \\
&$10^{-5}$&19&2573&309 &20&4153&275 &20&3094&204 &20&2632&176 &20&4194&278 \\ \hline
ES-D-ICCG&$10^{-3}$&20&2570&311 &20&4150&275 &20&3085&203 &20&2624&175 &20&4189&279 \\
&$10^{-4}$&20&2571&311 &20&4150&276 &20&3086&203 &20&2629&175 &20&4189&278 \\
&$10^{-5}$&19&2571&309 &20&4150&275 &20&3086&204 &20&2624&175 &20&4189&279 \\ \hline

\end{tabular}
\end{table*}

\begin{table*}
%\normalsize
\fontsize{9pt}{0.4cm}\selectfont

\centering
\caption{Numerical results (sequential solver, $\bbb$: random vector)}
\label{result-1-rand}
\begin{tabular}{lllllllllllllllll} \hline
 & & \multicolumn{3}{l}{Queen\_4147} & \multicolumn{3}{l}{Bump\_2911} & \multicolumn{3}{l}{G3\_circuit} & \multicolumn{3}{l}{Flan\_1565} & \multicolumn{3}{l}{Hook\_1498} \\	\cline{3-17}
Solver & $\theta$ & $\tilde{m}$ & \#Ite. & $T_{t}$ & $\tilde{m}$ & \#Ite. & $T_{t}$ & $\tilde{m}$ & \#Ite. & $T_{t}$ & $\tilde{m}$ & \#Ite. & $T_{t}$ & $\tilde{m}$ & \#Ite. & $T_{t}$ \\ \hline
ICCG & - & - & 3140 &  2776 & - &1544 & 564 & - & 926 & 46.1 & - & 3196 & 1010 & -  & 1613 & 282 \\ \hline
ES-SC-ICCG & $10^{-3}$ & 20 & 2546 & 2648 & 20 & 906 & 428 & 19 & 865 & 88.8 & 20 & 1013 & 372 & 20 & 554 & 127 \\
 & $10^{-4}$ & 19&2566& 2652 & 17& 921 & 422 & 10& 891 & 70.6 &19 & 1048 & 382 &13& 	672& 143 \\
 & $10^{-5}$ & 7 & 2783&	2626 & 5 & 1114 &	441 &	0 &-&-& 9 &1524& 518 & 5 &1076 & 208 \\ \hline
ES-D-ICCG  & $10^{-3}$ & 20 & 2542 & 2652 & 20 & 900 & 424 & 19 & 863 & 88.2 & 20 &	1011 & 372 &20 & 553 & 126 \\
 &  $10^{-4}$ & 19 & 2561 & 2654 & 17& 917 & 418 & 10& 890 & 70.2 &	19& 1048& 383 & 13	& 671 & 141 \\
&  $10^{-5}$ & 7 & 2779& 2630& 5& 1113 & 439 &0 &-& -& 9& 1523 & 517 &5& 1075 & 206 \\ \hline
\end{tabular}
\vspace{1\baselineskip}

\begin{tabular}{lllllllllllllllll} \hline
 & & \multicolumn{3}{l}{StocF-1465} & \multicolumn{3}{l}{Geo\_1438} & \multicolumn{3}{l}{Serena} & \multicolumn{3}{l}{thermal2}& \multicolumn{3}{l}{ecology2} \\ \cline{3-17}
Solver &$\theta$&$\tilde{m}$&\#Ite.&$T_{t}$&$\tilde{m}$&\#Ite.&$T_{t}$&$\tilde{m}$&\#Ite.&$T_{t}$&$\tilde{m}$&\#Ite.&$T_{t}$&$\tilde{m}$&\#Ite.&$T_{t}$ \\ \hline
ICCG&&-&55799 &4714 &-&441 &81.3 &-&299 &58.1 &-&2261 &141 &-&1902 &51.7 \\ \hline 
ES-SC-ICCG&$10^{-3}$&20&29693 &4001 &15&252 &54.7 &7&242 &49.1 &20&959 &101 &20&853 &52.5 \\
&$10^{-4}$&20&29693 &4011 &2&385 &71.9 &0&-&-&17&1020 &101 &16&933 &51.5 \\
&$10^{-5}$&20&29693 &4007 &0&-&-&0&-&-&4&1526 &112 &5&1268 &47.8 \\ \hline
ES-D-ICCG&$10^{-3}$&20&29600 &3990 &15&251 &55.1 &7&241 &49.6 &20&957 &100 &20&850 &52.2 \\ 
&$10^{-4}$&20&29600 &3993 &2&384 &72.3 &0&-&-&17&1019 &101 &16&930 &51.2  \\
&$10^{-5}$&20&29600 &4002 &0&-&-&0&-&-&4&1524 &111 &5&1267 &47.4 \\ \hline
\end{tabular}
\vspace{1\baselineskip}

\begin{tabular}{lllllllllllllllll} \hline
 & & \multicolumn{3}{l}{bone010} & \multicolumn{3}{l}{ldoor} & \multicolumn{3}{l}{audikw\_1} & \multicolumn{3}{l}{Emilia\_923}  &\multicolumn{3}{l}{boneS10} \\ \cline{3-17}
Solver &$\theta$&$\tilde{m}$&\#Ite.&$T_{t}$&$\tilde{m}$&\#Ite.&$T_{t}$&$\tilde{m}$&\#Ite.&$T_{t}$&$\tilde{m}$&\#Ite.&$T_{t}$&$\tilde{m}$&\#Ite.&$T_{t}$ \\ \hline
ICCG&&-&4189 &804 &-&2143 &293 &-&2420 &533 &-&459 &54 &-&8515 &1274 \\ \hline
ES-SC-ICCG&$10^{-3}$&20&996 &225 &20&1230 &208 &19&858 &214 &20&266 &39 &20&2733 &492 \\ 
&$10^{-4}$&17&1060 &234 &16&1259 &204 &8&1220 &284 &18&276 &40 &20&2733 &494  \\
&$10^{-5}$&13&1288 &277 &3&1649 &236 &4&1604 &364 &5&371 &46 &20&2733 &494  \\ \hline
ES-D-ICCG&$10^{-3}$&20&989 &221 &20&1227 &206 &19&861 &214 &20&267 &40 &20&2728 &490 \\  
&$10^{-4}$&17&1053 &231 &16&1256 &203 &8&1207 &280 &18&276 &40 &20&2728 &490 \\
&$10^{-5}$&13&1281 &273 &3&1648 &234 &4&1579 &356 &5&370 &47 &20&2728 &490 \\ \hline
\end{tabular}
\vspace{1\baselineskip}

\begin{tabular}{lllllllllllllllll} \hline
 & & \multicolumn{3}{l}{PFlow\_742} & \multicolumn{3}{l}{tmt\_sym}  & \multicolumn{3}{l}{apache2}  & \multicolumn{3}{l}{Fault\_639} & \multicolumn{3}{l}{parabolic\_fem} \\ \cline{3-17}
Solver &$\theta$&$\tilde{m}$&\#Ite.&$T_{t}$&$\tilde{m}$&\#Ite.&$T_{t}$&$\tilde{m}$&\#Ite.&$T_{t}$&$\tilde{m}$&\#Ite.&$T_{t}$&$\tilde{m}$&\#Ite.&$T_{t}$ \\ \hline
ICCG&&-&32971 &3311 &-&1256 &35.8 &-&770 &16.8 &-&2172 &176 &-&1208 &20.1 \\ \hline
ES-SC-ICCG&$10^{-3}$&20&14148 &1774 &20&562 &29.8 &20&342 &15.7 &20&1601 &162 &18&835 &27.6 \\ 
&$10^{-4}$&20&14148 &1772 &15&617 &29.5 &12&445 &16.4 &16&1629 &159 &9&889 &22.8 \\
&$10^{-5}$&20&14148 &1769 &3&1002 &33.9 &2&653 &16.7 &4&1899 &162 &0&-&- \\ \hline
ES-D-ICCG&$10^{-3}$&20&14042 &1758 &20&556 &29.8 &20&341 &15.7 &20&1595 &163 &18&833 &27.6 \\ 
&$10^{-4}$&20&14042 &1758 &15&614 &29.5 &12&444 &16.4 &16&1623 &159 &9&888 &22.6  \\
&$10^{-5}$&20&14042 &1754 &3&1000 &33.8 &2&653 &16.6 &4&1896 &162 &0&-&- \\
\end{tabular}
\vspace{1\baselineskip}

\begin{tabular}{lllllllllllllllll} \hline
 & & \multicolumn{3}{l}{bundle\_adj} & \multicolumn{3}{l}{af\_shell8}  & \multicolumn{3}{l}{af\_shell4} & \multicolumn{3}{l}{af\_shell3} & \multicolumn{3}{l}{af\_shell7} \\ \cline{3-17}
Solver &$\theta$&$\tilde{m}$&\#Ite.&$T_{t}$&$\tilde{m}$&\#Ite.&$T_{t}$&$\tilde{m}$&\#Ite.&$T_{t}$&$\tilde{m}$&\#Ite.&$T_{t}$&$\tilde{m}$&\#Ite.&$T_{t}$ \\ \hline
ICCG&&-&43578 &2325 &-&1038 &51.2 &-&1039 &51.7 &-&1039 &51.1 &-&1038 &51.6 \\ \hline 
ES-SC-ICCG&$10^{-3}$&20&11997 &846 &18&510 &32.8 &18&513 &33.5 &18&513 &33.2 &18&510 &33.2 \\ 
&$10^{-4}$&19&11394 &795 &9&606 &34.9 &9&602 &34.9 &9&602 &34.7 &9&606 &35.2 \\
&$10^{-5}$&19&11394 &795 &0&-&-&0&-&-&0&-&-&0&-&-\\ \hline
ES-D-ICCG&$10^{-3}$&20&10110 &711 &18&508 &33.0 &18&511 &33.4 &18&511 &33.1 &18&508 &33.1 \\ 
&$10^{-4}$&19&10030 &700 &9&605 &34.9 &9&601 &34.9 &9&601 &34.7 &9&605 &35.1 \\
&$10^{-5}$&19&10030 &698 &0&-&-&0&-&-&0&-&-&0&-&- \\ \hline
\end{tabular}
\vspace{1\baselineskip}

\begin{tabular}{lllllllllllllllll} \hline
 & & \multicolumn{3}{l}{inline\_1} & \multicolumn{3}{l}{af\_0\_k101}  & \multicolumn{3}{l}{af\_4\_k101}  & \multicolumn{3}{l}{af\_3\_k101}  & \multicolumn{3}{l}{af\_2\_k101} \\ \cline{3-17}
Solver &$\theta$&$\tilde{m}$&\#Ite.&$T_{t}$&$\tilde{m}$&\#Ite.&$T_{t}$&$\tilde{m}$&\#Ite.&$T_{t}$&$\tilde{m}$&\#Ite.&$T_{t}$&$\tilde{m}$&\#Ite.&$T_{t}$ \\ \hline
ICCG&&-&8464 &870 &-&12961 &641 &-&9974 &495 &-&8501 &420 &-&12970 &641 \\ \hline 
ES-SC-ICCG&$10^{-3}$&20&2686 &324 &20&5657 &372 &20&3582 &238 &20&2680 &177 &20&5413 &361 \\ 
&$10^{-4}$&20&2686 &324 &20&5657 &372 &20&3582 &238 &20&2680 &177 &20&5413 &360  \\
&$10^{-5}$&19&2686 &322 &20&5657 &372 &20&3582 &238 &20&2680 &177 &20&5413 &360  \\ \hline
ES-D-ICCG&$10^{-3}$&20&2683 &324 &20&5649 &376 &20&3575 &238 &20&2676 &177 &20&5409 &356 \\
&$10^{-4}$&20&2683 &324 &20&5649 &376 &20&3575 &238 &20&2676 &177 &20&5409 &356 \\
&$10^{-5}$&19&2684 &322 &20&5649 &376 &20&3575 &237 &20&2676 &177 &20&5409 &355 \\ \hline
\end{tabular}
\end{table*}

\subsection{Numerical results on the parallel solver}
In this subsection, we report the results of parallel (multithreaded) solver.
The parallelization of the CG solver is relatively straightforward.
But, the  IC preconditioning step that consists of forward and backward substitutions is not naturally parallelized.
While there are various parallel processing methods, we use one of simple but popular methods, namely block Jacobi IC preconditioning~\cite{saad} in the present research.
The parallelization of the subspace correction preconditioning and the deflation method is relatively easy.
The computationally dominant part for these methods is dense matrix vector multiplication, which can be straightforwardly parallelized.
Because $\tilde{m}$ is typically tiny, we sequentially process the solution part for the linear system with an $\tilde{m} \times \tilde{m} $ coefficient matrix $W^{\top}AW$ that is involved in the methods.

\begin{figure*}[tbp]
\centering
\includegraphics[clip, scale=0.9,  bb=50 90 560 400]{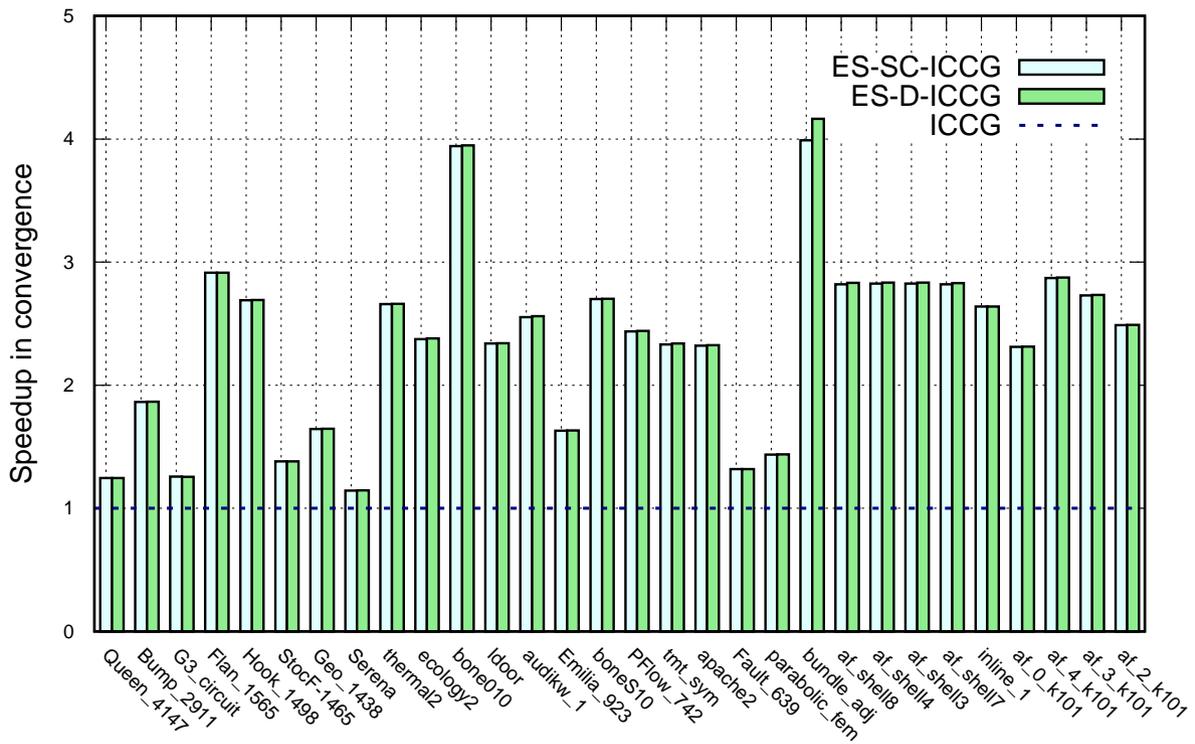}
\caption{Speedup in convergence of ES-SC-ICCG and ES-D-ICCG over ICCG (Parallel multithreaded solver, \bbb: random vector)}
\label{speedup-ite-para}
\end{figure*}

Tables \ref{result-p-ones} and \ref{result-p-rand} list the numerical results of the parallel ICCG solver and its variants with the proposed techniques, when a vector of ones and random vectors are used, respectively. From the viewpoint of convergence, the results on the parallel solver are similar to those of the sequential solver.
For all 60 test cases (30 datasets $\times$ 2 kinds of right-hand side vectors), convergence acceleration was attained by the proposed method. 
When a vector of ones was used, the convergence was more than twice as fast as the parallel ICCG solver for 27 out of 30 datasets.
Figure \ref{speedup-ite-para} shows the speedup in the convergence of the parallel solver based on the proposed technique against the parallel ICCG solver, when random vectors are used for the right-hand side vectors. In the test using random vectors,  the proposed method attains more than 2-fold speedup in convergence compared to the parallel ICCG solver for 21 out of 30 datasets.

\begin{figure*}[tbp]
\centering
\includegraphics[clip, scale=0.9,  bb=50 90 560 330]{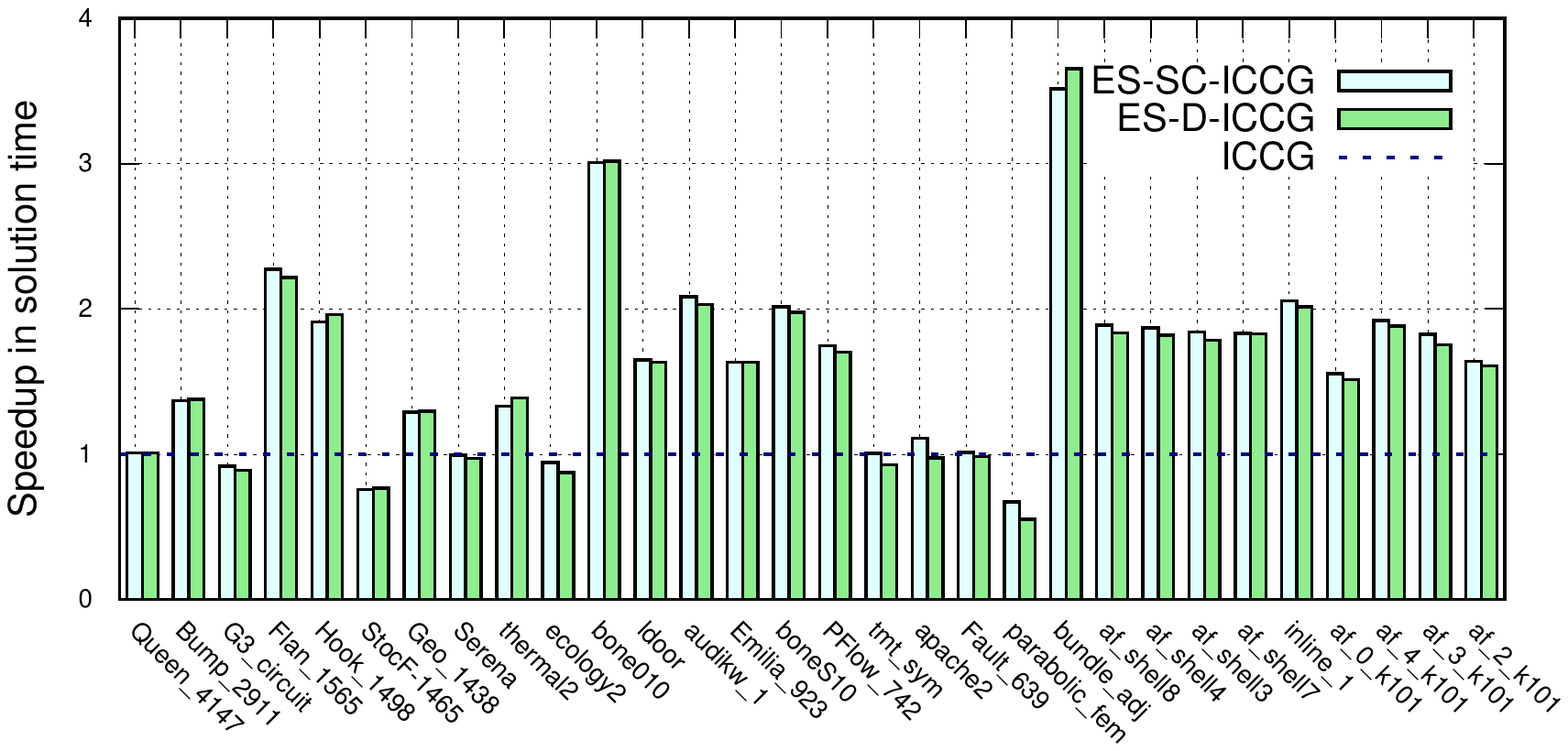}
\caption{Speedup in computational time of ES-SC-ICCG and ES-D-ICCG over pICCG (Parallel multithreaded solver, \bbb: random vector)}
\label{speedup-para}
\end{figure*}

Next, we examine the computational time. 
In the test using a vector of ones, the proposed method reduces the solution time for 28 out of 30 datasets.
For the bundle\_adj dataset, the parallel deflated ICCG solver based on our technique attains more than 4-fold speedup compared to the parallel ICCG solver.
The test using random vectors also indicates that our technique is effective to reduce the computational time for most of test datasets (25 out of 30).
In block Jacobi IC preconditioning, the computational cost for a PCG iteration is reduced as the number of threads increases.
Consequently, the impact of the additional cost for the convergence acceleration (subspace correction preconditioning or deflation) on the preconditioned solver is substantially enlarged in the parallel execution by many threads. In other words, the ratio of the iteration costs is enlarged from (\ref{ganma}).  Because we used a number of threads (=$\,$40) in our numerical tests, it becomes difficult to reduce the solution time compared with the sequential solver. However, Fig. \ref{speedup-para} indicates that our convergence acceleration technique accelerates the solution process for most of test problems.

%Table \ref{result-1-ones} lists the numerical results of the standard ICCG solver and its variants with the proposed convergence acceleration techniques, when a vector of ones is used for the right-hand side. ES-SC-ICCG denotes the solver of a CG solver with IC and subspace correction preconditioning based on proposed error vector sampling. The table show the average computational time in 50 solution steps. Table \ref{result-1-rand} shows the results when random vectors are used for the right-hand side. The table shows the average number of iterations and computational time in 50 solution steps. The numerical result indicates that both solvers based on the proposed method achieve convergence acceleration for all sixty test cases (30 datasets $\times$ 2 kinds of right-hand side vectors).  The convergence acceleration was significant for some of datasets. In the numerical tests using the vector of ones, the acceleration method attains more than 3-fold speedup in convergence for 16 our of 30 datasets.
%Even when we used random right-hand side vectors, the convergence rate was more than doubled for 20 out of 30 datasets.
%In this test, ES-SC-ICCG and ES-D-ICCG solvers obtained their best results for 12 out of 30 datasets when $\tilde{m}$ is equal to $m=20$.
%For these datasets, increase in the number of sample vectors, $m$ possibly improves the solver performance.
%Finally, Fig. \ref{speedup-seq} shows the speedup ratio of ES-SC-ICCG and ES-D-ICCG over ICCG, which confirms the effectiveness of our technique.

\subsection{Condition number estimation}
Figure \ref{eigen} implies that our technique based on error vector sampling can be a useful tool for the estimation of the smallest eigenvalue.
Because the estimation of the largest eigenvalue is relatively easy, the technique can be used for the estimation of the condition number of the coefficient matrix. 
Algorithm \ref{alg:PCG-con} shows the proposed procedure of PCG method with a condition number estimation.
The largest eigenvalue is estimated by the power method which is combined with the procedure of CG method.
The smallest eigenvalue is estimated by our technique.
We conducted numerical tests using five relatively small matrices downloaded from SuiteSparse matrix collection to examine our technique. Diagonal scaling is applied to the matrices before the tests. 
Table \ref{con1} shows the estimation for the largest eigenvalue $\lambda_{max}$, the smallest eigenvalue $\lambda_{min}$, and the condition number $\kappa$ compared with the results calculated by the LAPACK library.
Table \ref{con1} implies that the technique based on the error sampling is effective for the estimation of the condition number.

Because the number of sample vectors is much smaller than $n$ ($m \ll n$), the additional computational cost for the calculation of the smallest eigenvalue (Ritz value) is typically negligible compared with the iterative solution cost.
Although the power method requires an additional sparse matrix vector multiplication (SpMV) operation, it is combined with the SpMV for CG method. 
In this case, the matrix data transferred  from main memory are efficiently used for two vectors. Because the SpMV operation is typically memory bound, the impact of the power method on the iteration time is possibly limited.
Most of iterative solvers like CG solver usually uses a convergence criterion based on a (relative) residual norm.
If the estimation of the condition number of the coefficient matrix is given with the solution vector by the iterative solver, it can be a useful tool to evaluate the accuracy of the solution vector.
The proposed solver provides this function without large amount of additional computations.

\begin{algorithm}[t]
\caption{PCG method with condition number estimation}
\label{alg:PCG-con} 
\begin{algorithmic}[1]
\Require $A$, $\bmm{b}$, $M$, $\bmm{x}_{0}$, $\varepsilon$, $m$
\State $\bmm{r}_0 =  \bmm{b} - A \bmm{x}_0 $
%\State $\bmm{p}_0 = M^{-1} \bmm{r}_0$
\State $\bmm{v}_{0} \leftarrow$ Initialization (a unit random vector etc.) 
%\For{$i = 0, 1, 2, \ldots $ \textbf{until} $ \| \bmm{r}_i \|_2 \le \varepsilon \| \bmm{b} \|_2$}
\For{$i = 1, 2, \ldots $}
\State $\bmm{z}_{i-1} = M^{-1} \bmm{r}_{i-1}$
\State $\rho_{i-1}=(\rrr_{i-1}, \zzz_{i-1})$
\If{i=1}
\State $\ppp_{1}=\zzz_{0}$
\Else
\State $\beta_{i-1}=\rho_{i-1}/\rho_{i-2}$
\State $\ppp_{i}=\zzz_{i-1}+\beta_{i-1}\ppp_{i-1}$
\EndIf
\State $(\bmm{q}_{i} \ \bmm{v}_{i}) = A (\bmm{p}_{i} \ \bmm{v}_{i-1}) $ \ \ \ // (SpMV)
\State $\bmm{v}_{i} = \bmm{v}_{i} / \| \bmm{v}_{i} \|_2$
\State $\alpha_{i}=\rho_{i-1}/(\ppp_{i}, \qqq_{i})$
\State $\bmm{x}_{i} = \bmm{x}_{i-1} + \alpha_i \bmm{p}_{i}$
\State $\bmm{r}_{i} = \bmm{r}_{i-1} - \alpha_i \bmm{q}_{i}$
\If{$\| \bmm{r}_{i} \|_2 \le \varepsilon \| \bmm{b} \|_2$} 
\State break
\EndIf
%
%\If{Sampling condition is satisfied} 
%\State $\tilde{\xxx}^{(s)}=\bmm{x}_{i}, \ s \in \{1, 2, \ldots, m \}$ 
%\EndIf
%\State $(\bmm{q}_{i} \ \bmm{v}_{i}) = A (\bmm{p}_{i-1} \ \bmm{v}_{i-1}) $ \ \ \ // (SpMV)
%\State $\bmm{v}_{i} = \bmm{v}_{i} / \| \bmm{v}_{i} \|$
%\State $\displaystyle \alpha_i = \frac{(M^{-1}\bmm{r}_{i-1}, \bmm{r}_{i-1})}{(\bmm{p}_{i-1}, \bmm{q}_{i})}$
%\State $\bmm{x}_{i} = \bmm{x}_{i-1} + \alpha_i \bmm{p}_{i-1}$
%\State $\bmm{r}_{i} = \bmm{r}_{i-1} - \alpha_i \bmm{q}_{i}$
%\If{$\| \bmm{r}_{i} \|_2 \le \varepsilon \| \bmm{b} \|_2$} 
%\State break
%\EndIf
\If{Sampling condition is satisfied} 
\State $\tilde{\xxx}^{(s)}=\bmm{x}_{i}, \ s \in \{1, 2, \ldots, m \}$ 
\EndIf
%\State $\displaystyle \beta_i = - \frac{(M^{-1}\bmm{r}_{i}, \bmm{r}_{i})}{(M^{-1}\bmm{r}_{i-1}, \bmm{r}_{i-1})}$
%\State $\bmm{p}_{i} =  M^{-1} \bmm{r}_{i} + \beta_i \bmm{p}_{i-1}$
%
%\If{$\| \bmm{r}_{i+1} \|_2 \le \varepsilon \| \bmm{b} \|_2$} 
%break
%\EndIf
\EndFor
\State $\tilde{E}=(\bmm{x}_{i}-\tilde{\xxx}^{(1)} \, \bmm{x}_i-\tilde{\xxx}^{(2)} \,  \ldots \, \bmm{x}_i-\tilde{\xxx}^{(m)} )$
\State Apply Gram-Schmidt method to $\tilde{E}$ and obtain $E$
\State Solve an eigenvalue problem:  $E^{\top} A E \ttt = \lambda \ttt$ and obtain the smallest Ritz value $\lambda_{min}$
\State $\lambda_{max}=(\bmm{v}_{i}, A\bmm{v}_{i})$
\State $\kappa=\lambda_{max}/\lambda_{min}$
\Ensure $\bmm{x}_{i}$, $\kappa$
\end{algorithmic}
\end{algorithm}

\section{Conclusion}
In this paper, we introduce an algebraic auxiliary matrix construction method that can be used for the subspace correction preconditioning and the deflation method. We focus on a problem in which a sequence of linear systems with an identical coefficient matrix are solved.
In our method, we sample the approximate solution vectors in the first iterative solution step, and calculate the error vectors corresponding to the sample vectors after the solution step is completed. 
Then, we perform the Rayleigh-Ritz method using a subspace spanned by these error vectors to identify (approximate) eigenvectors associated with small eigenvalues.
Finally, the auxiliary matrix is constructed by the Ritz vectors associated with small Ritz values.
We also present a cost model of the subspace preconditioning and the deflation method.
Numerical tests using 30 coefficient matrices were conducted to verify our technique.
The test results confirm that the proposed convergence acceleration technique efficiently reduces both the number of iterations for the convergence and the solution time of the serial and parallel preconditioned CG solvers. 
Moreover, additional numerical tests indicate that the proposed technique can be used for the condition number estimation. 

Currently, we examine the effectiveness of the technique for a linear system having an unsymmetric coefficient matrix. Because the preliminary results show its effectiveness, we will report it in future. We are also investigating application of the technique to other problems. Especially, we examine the effectiveness in parallel-in-time (PinT) simulations, which often involves the solution process of multiple linear systems of coefficient matrices having a common property. We are also interested in the combination of our technique with preconditioning techniques suitable for GPU computing. In future, we will examine our technique in various situations of computational science or engineering problems.

%In this paper, we proposed a mapping operator generation method for subspace correction to accelerate transient and/or nonlinear eddy-current analyses. In our method, the operator is automatically generated using the information obtained in the preceding time or nonlinear iteration steps. 
%Transient nonlinear analyses with the transformer and the box shield models demonstrated that the proposed method attained a significant reduction in the number of iterations and computational time, respectively.
%Although the numerical result is not shown here because of a lack of space, we also confirmed the effectiveness of the proposed method in transient linear analyses.
% 
%To the best of our knowledge, the work by Gosselet et al.~\cite{Gosselet} is the most
%closely related to ours. However, in their method, all residual vectors in previous solution steps are preserved, which requires huge additional memory space.
%Consequently, our method has significantly different memory requirements from Gosselet's method, which is important in the application to practical analyses.
%
%In future work, we will examine our method for more numerical test models with various time-step sizes. Moreover, we will investigate the effect of the method when it is used with various nonlinear and linear iteration methods, for example, fixed point iterations and a multigrid method.
%

\begin{table*}
%\normalsize
\centering
\caption{Solver performance using sampling method B (sequential solver, $\bbb=(1, 1, \ldots, 1)^\top$)}
\label{sampling B}
\begin{tabular}{llllllll} \hline
& & \multicolumn{3}{l}{Flan\_1565} & \multicolumn{3}{l}{Hook\_1498} \\	\cline{3-8}
Solver & $\theta$ & $\tilde{m}$ & \#Ite. & $T_{t}$ & $\tilde{m}$ & \#Ite. & $T_{t}$ \\ \hline
ES-SC-ICCG&$10^{-3}$&20&1584 & 586 &15&1075 &233  \\ 
&$10^{-4}$&15&1927 &690&7&1157 &229 \\
&$10^{-5}$&7&2094 & 706 &4&1208 &230  \\ \hline
ES-D-ICCG&$10^{-3}$ & 20 & 1579 & 585 & 15 & 1072 & 232\\ 
&$10^{-4}$&15&1925&687 &7& 1156&229 \\
&$10^{-5}$&7&2093 &704 &5& 1207 &229\\ \hline
\end{tabular}
\end{table*}

\begin{table*}
%\normalsize
\fontsize{9pt}{0.4cm}\selectfont

\centering
\caption{Numerical results (parallel solver, $\bbb=(1, 1, \ldots, 1)^\top$)}
\label{result-p-ones}
\begin{tabular}{lllllllllllllllll} \hline
 & & \multicolumn{3}{l}{Queen\_4147} & \multicolumn{3}{l}{Bump\_2911} & \multicolumn{3}{l}{G3\_circuit} & \multicolumn{3}{l}{Flan\_1565} & \multicolumn{3}{l}{Hook\_1498} \\	\cline{3-17}
Solver & $\theta$ & $\tilde{m}$ & \#Ite. & $T_{t}$ & $\tilde{m}$ & \#Ite. & $T_{t}$ & $\tilde{m}$ & \#Ite. & $T_{t}$ & $\tilde{m}$ & \#Ite. & $T_{t}$ & $\tilde{m}$ & \#Ite. & $T_{t}$ \\ \hline
ICCG&&-&4663 &215 &-&3455 &73.4 &-&1461 &5.40 &-&4911 &86.6 &-&2312 &23.4 \\ \hline
ES-SC-ICCG&$10^{-3}$&20&1532 &89 &20&1062 &30.8 &20&468 &3.81 &20&1504 &31.3 &19.0 &808 &11.6 \\
&$10^{-4}$&20&1532&88 &17&1814&51.1 &13&1067&7.22 &18&1777 &37.3 &13.0 &1036 &13.9 \\
&$10^{-5}$&6&4258&218 &2&2977&68.1 &2&1392&5.67 &9&2415 &47.2 &5.0 &1562 &18.8 \\ \hline
ES-D-ICCG&$10^{-3}$&20&1530&91 &20&1062&31.3 &20&468&3.90 &20&1502&33.0 &19.0 &807&12.1 \\
&$10^{-4}$&20&1530&93 &17&1811&51.4 &13&1066&7.08 &18&1776&38.2 &13.0 &1035&14.3 \\
&$10^{-5}$&6&4252 &216 &2&2977 &68.7 &2&1392 &5.73 &9&2409 &46.6 &5.0 &1561 &18.9 \\ \hline
\end{tabular}
\vspace{1\baselineskip}

\begin{tabular}{lllllllllllllllll} \hline
 & & \multicolumn{3}{l}{StocF-1465} & \multicolumn{3}{l}{Geo\_1438} & \multicolumn{3}{l}{Serena} & \multicolumn{3}{l}{thermal2}& \multicolumn{3}{l}{ecology2} \\ \cline{3-17}
Solver &$\theta$&$\tilde{m}$&\#Ite.&$T_{t}$&$\tilde{m}$&\#Ite.&$T_{t}$&$\tilde{m}$&\#Ite.&$T_{t}$&$\tilde{m}$&\#Ite.&$T_{t}$&$\tilde{m}$&\#Ite.&$T_{t}$ \\ \hline
ICCG&&-&66348&329 &-&904&9.62 &-&628&6.77 &-&3583&12.0 &-&2131&3.29 \\ \hline
ES-SC-ICCG&$10^{-3}$&20&16453&157 &14&549&7.26 &8&546&6.36 &20&1128&8.1 &20&885&4.40 \\
&$10^{-4}$&20&16453&154 &2&779&8.24 &0&-&-&17&1555&10.7 &15&1039&4.19 \\
&$10^{-5}$&20&16453&157 &0&-&-&0&-&-&4&2506&10.4 &4&1656&3.95 \\ \hline
ES-D-ICCG&$10^{-3}$&20&16452&148 &14&548&7.27 &8&545&6.39 &20&1126&8.0 &20&882&4.38 \\
&$10^{-4}$&20&16452&147 &2&778&8.40 &0&-&-&17&1554&10.2 &15&1037&4.33 \\
&$10^{-5}$&20&16452&150 &0&-&-&0&-&-&4&2504&10.8 &4&1654&4.05 \\ \hline
\end{tabular}
\vspace{1\baselineskip}

\begin{tabular}{lllllllllllllllll} \hline
 & & \multicolumn{3}{l}{bone010} & \multicolumn{3}{l}{ldoor} & \multicolumn{3}{l}{audikw\_1} & \multicolumn{3}{l}{Emilia\_923}  &\multicolumn{3}{l}{boneS10} \\ \cline{3-17}
Solver &$\theta$&$\tilde{m}$&\#Ite.&$T_{t}$&$\tilde{m}$&\#Ite.&$T_{t}$&$\tilde{m}$&\#Ite.&$T_{t}$&$\tilde{m}$&\#Ite.&$T_{t}$&$\tilde{m}$&\#Ite.&$T_{t}$ \\ \hline
ICCG&&-&7838 &76.9 &-&5227 &38.0 &-&2635 &30.0 &-&6542 &42.1 &-&14690 &119 \\ \hline
ES-SC-ICCG&$10^{-3}$&20&2141 &29.6 &20&1503 &14.8 &20&816 &11.4 &20&1893 &17.3 &20&5166 &55 \\
&$10^{-4}$&18&2207&28.6 &18&2199&20.7 &7&1549&19.1 &19&2991&27.5 &20&5164 &56 \\
&$10^{-5}$&12&2925&35.7 &3&4040&29.6 &3&1798&21.2 &7&4829&36.0 &19&5446 &58 \\ \hline
ES-D-ICCG&$10^{-3}$&20&2138&28.9 &20&1504&14.8 &20&813&11.6 &20&1895&17.8 &20&5162&57 \\
&$10^{-4}$&18&2203&28.9 &18&2196&20.9 &7&1545&19.3 &19&2989&29.1 &20&5161&56 \\
&$10^{-5}$&12&2921 &36.3 &3&4036 &30.8 &3&1796 &21.9 &7&4823 &37.1 &19&5446 &59 \\ \hline
\end{tabular}
\vspace{1\baselineskip}

\begin{tabular}{lllllllllllllllll} \hline
 & & \multicolumn{3}{l}{PFlow\_742} & \multicolumn{3}{l}{tmt\_sym}  & \multicolumn{3}{l}{apache2}  & \multicolumn{3}{l}{Fault\_639} & \multicolumn{3}{l}{parabolic\_fem} \\ \cline{3-17}
Solver &$\theta$&$\tilde{m}$&\#Ite.&$T_{t}$&$\tilde{m}$&\#Ite.&$T_{t}$&$\tilde{m}$&\#Ite.&$T_{t}$&$\tilde{m}$&\#Ite.&$T_{t}$&$\tilde{m}$&\#Ite.&$T_{t}$ \\ \hline
ICCG&&-&37485&214 &-&1576&2.26 &-&1056&1.31 &-&5083&26.2 &-&2125&1.58 \\ \hline
ES-SC-ICCG&$10^{-3}$&20&11633&95 &20&638&2.40 &19&408&1.51 &20&1496&9.6 &19&1419&3.50 \\ 
&$10^{-4}$&20&11633&95 &14&777&2.33 &12&494&1.27 &18&3075&19.4 &8&1326&1.89 \\
&$10^{-5}$&20&11633&99 &3&1259&2.09 &3&816&1.31 &2&4735&22.6 &0&-&- \\ \hline
ES-D-ICCG&$10^{-3}$&20&11617&100 &20&636&2.44 &19&408&1.53 &20&1495&9.6 &19&1417&4.03 \\
&$10^{-4}$&20&11617&97 &14&776&2.35 &12&494&1.39 &18&3074&18.8 &8&1325&2.27 \\
&$10^{-5}$&20&11617&102 &3&1257&2.30 &3&815&1.43 &2&4739&22.2 &0&-&- \\ \hline
\end{tabular}
\vspace{1\baselineskip}

\begin{tabular}{lllllllllllllllll} \hline
 & & \multicolumn{3}{l}{bundle\_adj} & \multicolumn{3}{l}{af\_shell8}  & \multicolumn{3}{l}{af\_shell4} & \multicolumn{3}{l}{af\_shell3} & \multicolumn{3}{l}{af\_shell7} \\ \cline{3-17}
Solver &$\theta$&$\tilde{m}$&\#Ite.&$T_{t}$&$\tilde{m}$&\#Ite.&$T_{t}$&$\tilde{m}$&\#Ite.&$T_{t}$&$\tilde{m}$&\#Ite.&$T_{t}$&$\tilde{m}$&\#Ite.&$T_{t}$ \\ \hline
ICCG&&-&64356 &797 &-&1575 &5.13 &-&1575 &5.31 &-&1575 &5.07 &-&1575 &4.91  \\ \hline
ES-SC-ICCG&$10^{-3}$&20&14407 &208 &20&537 &2.58 &20&539 &2.59 &20&539 &2.62 &20&537 &2.56 \\
&$10^{-4}$&20&14407 &207 &11&764 &3.08 &11&764 &3.10 &11&764 &3.17 &11&764 &3.07 \\
&$10^{-5}$&19&14104 &200 &0&-&-&0&-&-&0&-&-&0&-&- \\ \hline
ES-D-ICCG&$10^{-3}$&20&13547&196 &20&537&2.64 &20&537&2.76 &20&537&2.65 &20&537&2.59 \\
&$10^{-4}$&20&13547&196 &11&764&3.25 &11&763&3.33 &11&763&3.24 &11&764&3.20 \\
&$10^{-5}$&19&13611 &194 &0&-&-&0&-&-&0&-&-&0&-&- \\ \hline
\end{tabular}
\vspace{1\baselineskip}

\begin{tabular}{lllllllllllllllll} \hline
 & & \multicolumn{3}{l}{inline\_1} & \multicolumn{3}{l}{af\_0\_k101}  & \multicolumn{3}{l}{af\_4\_k101}  & \multicolumn{3}{l}{af\_3\_k101}  & \multicolumn{3}{l}{af\_2\_k101} \\ \cline{3-17}
Solver &$\theta$&$\tilde{m}$&\#Ite.&$T_{t}$&$\tilde{m}$&\#Ite.&$T_{t}$&$\tilde{m}$&\#Ite.&$T_{t}$&$\tilde{m}$&\#Ite.&$T_{t}$&$\tilde{m}$&\#Ite.&$T_{t}$ \\ \hline
ICCG&&-&23064&115 &-&16157&48.6 &-&12458&45.5 &-&10595&34.9 &-&16249&57.6  \\ \hline
ES-SC-ICCG&$10^{-3}$&20&6393&44 &20&5026&24.7 &20&4230&20.7 &20&3567&17.1 &20&4924&24.6 \\
&$10^{-4}$&20&6393&43 &20&5026&23.9 &20&4230&19.7 &20&3567&16.9 &20&4924&23.7 \\
&$10^{-5}$&18&8969&60 &20&5026&24.1 &20&4230&19.5 &20&3567&16.8 &20&4924&22.8  \\ \hline
ES-D-ICCG&$10^{-3}$&20&6390&46 &20&5018&24.8 &20&4237&21.6 &20&3574&18.0 &20&4920&24.9 \\ 
&$10^{-4}$&20&6390&46 &20&5018&24.7 &20&4237&20.8 &20&3574&17.6 &20&4920&24.5 \\
&$10^{-5}$&18&8964&62 &20&5018&25.2 &20&4237&20.4 &20&3574&17.5 &20&4920&24.4  \\ \hline
\end{tabular}
\end{table*}

\begin{table*}
%\normalsize
\fontsize{9pt}{0.4cm}\selectfont

\centering
\caption{Numerical results (parallel solver, $\bbb$: random vector)}
\label{result-p-rand}
\begin{tabular}{lllllllllllllllll} \hline
 & & \multicolumn{3}{l}{Queen\_4147} & \multicolumn{3}{l}{Bump\_2911} & \multicolumn{3}{l}{G3\_circuit} & \multicolumn{3}{l}{Flan\_1565} & \multicolumn{3}{l}{Hook\_1498} \\	\cline{3-17}
Solver & $\theta$ & $\tilde{m}$ & \#Ite. & $T_{t}$ & $\tilde{m}$ & \#Ite. & $T_{t}$ & $\tilde{m}$ & \#Ite. & $T_{t}$ & $\tilde{m}$ & \#Ite. & $T_{t}$ & $\tilde{m}$ & \#Ite. & $T_{t}$ \\ \hline
ICCG&&-&4684 &215 &-&3437 &73.9 &-&1455 &5.27 &-&4906 &83.6 &-&2309 &24.3 \\ \hline
ES-SC-ICCG&$10^{-3}$&20&3762 &217 &20&1844 &53.9 &20&1157 &9.54 &20&1684 &36.8 &19.0 &858 &12.7 \\
&$10^{-4}$&19&3760 &217 &17&1929 &54.6 &13&1203 &8.15 &18&1768 &37.8 &13.0 &1027 &13.6 \\
&$10^{-5}$&6&4166 &212 &2&2967 &67.8 &2&1388 &5.74 &9&2411 &47.3 &5.0 &1557 &18.3 \\ \hline
ES-D-ICCG&$10^{-3}$&20&3761 &220 &20&1842 &53.5 &20&1159 &9.71 &20&1684 &37.7 &19&857 &12.4 \\ 
&$10^{-4}$&19&3758 &218 &17&1927 &53.8 &13&1202 &8.24 &18&1766 &38.5 &13&1026 &13.9 \\
&$10^{-5}$&6&4164 &212 &2&2964 &67.0 &2&1388 &5.94 &9&2410 &48.2 &5.0 &1556 &18.6 \\ \hline
\end{tabular}
\vspace{1\baselineskip}

\begin{tabular}{lllllllllllllllll} \hline
 & & \multicolumn{3}{l}{StocF-1465} & \multicolumn{3}{l}{Geo\_1438} & \multicolumn{3}{l}{Serena} & \multicolumn{3}{l}{thermal2}& \multicolumn{3}{l}{ecology2} \\ \cline{3-17}
Solver &$\theta$&$\tilde{m}$&\#Ite.&$T_{t}$&$\tilde{m}$&\#Ite.&$T_{t}$&$\tilde{m}$&\#Ite.&$T_{t}$&$\tilde{m}$&\#Ite.&$T_{t}$&$\tilde{m}$&\#Ite.&$T_{t}$ \\ \hline
ICCG&&-&51167 &258 &-&901 &9.37 &20&625 &6.49 &-&3569 &12.5 &-&2225 &3.48 \\ \hline
ES-SC-ICCG&$10^{-3}$&20&37038 &341 &14&548 &7.26 &8&546 &6.53 &20&1343 &9.4 &20&937 &4.67 \\
&$10^{-4}$&20&37038 &341 &2&774 &8.49 &0&-&-&17&1597 &10.4 &17&1045 &4.74 \\
&$10^{-5}$&20&37038 &342 &0&-&-&0&-&-&4&2466 &10.6 &5&1471 &3.69 \\ \hline
ES-D-ICCG&$10^{-3}$&20&37036 &342 &14&547 &7.22 &8&545 &6.69 &20&1342 &9.1 &20&935 &4.55 \\
&$10^{-4}$&20&37036 &340 &2&773 &8.57 &0&-&-&17&1596 &10.2 &17&1044 &4.64 \\
&$10^{-5}$&20&37036 &337 &0&-&-&0&-&-&4&2465 &10.8 &5&1469 &3.99 \\ \hline
\end{tabular}
\vspace{1\baselineskip}

\begin{tabular}{lllllllllllllllll} \hline
 & & \multicolumn{3}{l}{bone010} & \multicolumn{3}{l}{ldoor} & \multicolumn{3}{l}{audikw\_1} & \multicolumn{3}{l}{Emilia\_923}  &\multicolumn{3}{l}{boneS10} \\ \cline{3-17}
Solver &$\theta$&$\tilde{m}$&\#Ite.&$T_{t}$&$\tilde{m}$&\#Ite.&$T_{t}$&$\tilde{m}$&\#Ite.&$T_{t}$&$\tilde{m}$&\#Ite.&$T_{t}$&$\tilde{m}$&\#Ite.&$T_{t}$ \\ \hline
ICCG&&&8190 &84.1 &-&5198 &36.2 &-&2633 &30.2 &-&6507 &43.0 &-&14637 &119  \\ \hline
ES-SC-ICCG&$10^{-3}$&20&2077 &28.3 &20&2222 &22.0 &20&1031 &14.5 &20&3989 &35.4 &20&5422 &60 \\
&$10^{-4}$&19&2100 &28.0 &18&2320 &21.9 &7&1541 &19.4 &19&4027 &36.7 &20&5422 &59 \\
&$10^{-5}$&12&2883 &35.6 &3&4019 &30.1 &3&1791 &21.7 &7&4798 &36.8 &19&5432 &59  \\ \hline
ES-D-ICCG&$10^{-3}$&20&2074 &27.9 &20&2220 &22.6 &20&1028 &14.9 &20&3986 &37.7 &20&5419 &61 \\ 
&$10^{-4}$&19&2096 &28.4 &18&2319 &22.1 &7&1536 &19.9 &19&4026 &37.9 &20&5419 &62 \\
&$10^{-5}$&12&2870 &37.1 &3&4016 &31.2 &3&1789 &21.8 &7&4797 &36.9 &19&5428 &60  \\ \hline
\end{tabular}
\vspace{1\baselineskip}

\begin{tabular}{lllllllllllllllll} \hline
 & & \multicolumn{3}{l}{PFlow\_742} & \multicolumn{3}{l}{tmt\_sym}  & \multicolumn{3}{l}{apache2}  & \multicolumn{3}{l}{Fault\_639} & \multicolumn{3}{l}{parabolic\_fem} \\ \cline{3-17}
Solver &$\theta$&$\tilde{m}$&\#Ite.&$T_{t}$&$\tilde{m}$&\#Ite.&$T_{t}$&$\tilde{m}$&\#Ite.&$T_{t}$&$\tilde{m}$&\#Ite.&$T_{t}$&$\tilde{m}$&\#Ite.&$T_{t}$ \\ \hline
ICCG&&-&37486 &221 &-&1569 &2.13 &-&1055 &1.34 &-&5047 &22.7 &-&2583 &1.89 \\ \hline
ES-SC-ICCG&$10^{-3}$&20&15380 &127 &20&673 &2.53 &19&454 &1.68 &20&3828 &24.9 &19&1798 &4.33 \\
&$10^{-4}$&20&15380 &127 &14&784 &2.51 &12&499 &1.37 &18&3872 &24.3 &9&1885 &2.80 \\
&$10^{-5}$&20&15380 &130 &3&1256 &2.12 &3&810 &1.21 &2&4710 &22.4 &0&-&-  \\ \hline
ES-D-ICCG&$10^{-3}$&20&15354 &131 &20&671 &2.55 &19&454 &1.64 &20&3830 &24.9 &19&1797 &4.84 \\ 
&$10^{-4}$&20&15354 &130 &14&782 &2.53 &12&498 &1.51 &18&3869 &24.6 &9&1885 &3.42 \\
&$10^{-5}$&20&15354 &133 &3&1255 &2.30 &3&809 &1.38 &2&4708 &23.0 &0&-&-  \\ \hline
\end{tabular}
\vspace{1\baselineskip}

\begin{tabular}{lllllllllllllllll} \hline
 & & \multicolumn{3}{l}{bundle\_adj} & \multicolumn{3}{l}{af\_shell8}  & \multicolumn{3}{l}{af\_shell4} & \multicolumn{3}{l}{af\_shell3} & \multicolumn{3}{l}{af\_shell7} \\ \cline{3-17}
Solver &$\theta$&$\tilde{m}$&\#Ite.&$T_{t}$&$\tilde{m}$&\#Ite.&$T_{t}$&$\tilde{m}$&\#Ite.&$T_{t}$&$\tilde{m}$&\#Ite.&$T_{t}$&$\tilde{m}$&\#Ite.&$T_{t}$ \\ \hline
ICCG&&-&55336 &701 &-&1572 &5.07 &-&1572 &5.06 &-&1572 &4.95 &-&1572 &5.04 \\ \hline 
ES-SC-ICCG&$10^{-3}$&20&13873 &199 &20&557 &2.69 &20&556 &2.71 &20&556 &2.68 &20&557 &2.75 \\
&$10^{-4}$&20&13873 &200 &11&760 &3.05 &11&760 &3.04 &11&760 &3.04 &11&760 &3.03 \\
&$10^{-5}$&19&13947 &200 &0&-&-&0&-&-&0&-&-&0&-&-  \\ \hline
ES-D-ICCG&$10^{-3}$&20&13287 &192 &20&555 &2.76 &20&555 &2.78 &20&555 &2.77 &20&555 &2.75 \\
&$10^{-4}$&20&13287 &193 &11&759 &3.23 &11&759 &3.21 &11&759 &3.19 &11&759 &3.22 \\
&$10^{-5}$&19&13796 &199 &0&-&-&0&-&-&0&-&-&0&-&-  \\ \hline
\end{tabular}
\vspace{1\baselineskip}

\begin{tabular}{lllllllllllllllll} \hline
 & & \multicolumn{3}{l}{inline\_1} & \multicolumn{3}{l}{af\_0\_k101}  & \multicolumn{3}{l}{af\_4\_k101}  & \multicolumn{3}{l}{af\_3\_k101}  & \multicolumn{3}{l}{af\_2\_k101} \\ \cline{3-17}
Solver &$\theta$&$\tilde{m}$&\#Ite.&$T_{t}$&$\tilde{m}$&\#Ite.&$T_{t}$&$\tilde{m}$&\#Ite.&$T_{t}$&$\tilde{m}$&\#Ite.&$T_{t}$&$\tilde{m}$&\#Ite.&$T_{t}$ \\ \hline
ICCG&&-&23054 &124 &-&16121 &51.5 &-&12425 &39.9 &-&10584 &33.7 &-&16237 &50.9 \\ \hline
ES-SC-ICCG&$10^{-3}$&20&8738 &60 &20&6977 &33.1 &20&4327 &20.8 &20&3877 &19.6 &20&6526 &31.0 \\ 
&$10^{-4}$&20&8738 &61 &20&6977 &33.1 &20&4327 &20.9 &20&3877 &18.8 &20&6526 &31.1 \\
&$10^{-5}$&18&9090 &62 &20&6977 &33.5 &20&4327 &20.5 &20&3877 &18.5 &20&6526 &31.3 \\ \hline
ES-D-ICCG&$10^{-3}$&20&8732 &61 &20&6966 &34.0 &20&4322 &21.2 &20&3873 &19.2 &20&6523 &32.1 \\ 
&$10^{-4}$&20&8732 &62 &20&6966 &34.6 &20&4322 &21.2 &20&3873 &19.3 &20&6523 &31.7 \\
&$10^{-5}$&18&9086 &65 &20&6966 &34.5 &20&4322 &21.3 &20&3873 &19.3 &20&6523 &32.2 \\ \hline
\end{tabular}
\end{table*}

\begin{table*}[tbp]
%      \small
	\centering
%	\vspace{0.75\baselineskip}
	\caption{Condition number estimation based on error vector sampling}
	\label{con1}
%	\begin{tabular}{|c|c|c|c|c|}
	\begin{tabular}{lllllll}
	\hline
	 & \multicolumn{3}{l}{Estimation} & \multicolumn{3}{l}{LAPACK} \\ \cline{2-7} 
Dataset & $\lambda_{max}$ & $\lambda_{min}$ & $\kappa$ & $\lambda_{max}$ & $\lambda_{min}$ & $\kappa$   \\ 
	\hline
bcsstk07 & 2.88 & 9.11$\times 10^{-5}$ & 3.16$\times 10^{4}$ & 2.90 & 9.11$\times 10^{-5}$ & 3.18$\times 10^{4}$ \\
msc01440 & 3.62& 3.08$\times 10^{-4}$ & 1.18$\times 10^{4}$ & 3.62 & 2.86$\times 10^{-4}$ & 1.27$\times 10^{4}$ \\
ex33 & 3.93 & 2.86$\times 10^{-11}$ & 1.37$\times 10^{11}$ & 3.93 & 2.59$\times 10^{11}$ & 1.52$\times 10^{11}$ \\
494\_bus & 1.99 & 2.55$\times 10^{-5}$ & 7.81$\times 10^{4}$ & 2.00 & 2.53$\times 10^{-5}$ & 7.90$\times 10^{4}$ \\
bcsstk06 & 2.89 &9.21$\times 10^{-5}$ & 3.14$\times 10^{4}$ & 2.90 & 9.11$\times 10^{-5}$ & 3.18$\times 10^{4}$ \\ \hline
\end{tabular}
\end{table*}

\appendix
\section{Selection method for approximation vectors\label{A}}
Algorithm \ref{alg1} shows the sampling method A for the approximate solution vector~\cite{magn}.
In the algorithm, $i$ is the iteration count, and $m$ is the number of sample vectors.
The parameter $l_{max}$ is set to satisfy $m^{l_{max}}>N_{max}$, where $N_{max}$ is the preset maximum iteration count of the solver.

\begin{algorithm}                      
\caption{Selection of approximate solution vectors}         
\label{alg1}                          
\begin{algorithmic}                  
\State $h=1$
\For{$i=1, 2, \cdots$}
\State {\bf Solver part}
\State Convergence check
\If  {$(\mbox{mod}(i, h)==0)$} 
\State  $i_{t}=\sum_{l=0}^{l_{max}} (-1)^{l} \lfloor (i-1)/m^{l} \rfloor $
\State  $s=\mbox{mod}(i_{t}, m)+1$
\State  $\tilde{\xxx}^{(s)}=\tilde{\xxx}_{i}$
\If {$(i==h*m)$} 
\State $h=h*2$ 
\EndIf
\EndIf
\EndFor
\end{algorithmic}
\end{algorithm}

%\begin{algorithm}                      
%\caption{Selection of approximate solution vectors}         
%\label{alg1}                          
%\begin{algorithmic}                  
%\State $h=1$
%\For{$i_{te}=1, 2, \cdots$}
%\State {\bf Solver part}
%\State Convergence check
%\If  {$(\mbox{mod}(i_{te}, h)==0)$} 
%\State  $i_{t}=\sum_{l=0}^{l_{max}} (-1)^{l} \lfloor (i_{te}-1)/m^{l} \rfloor $
%\State  $j=\mbox{mod}(i_{t}, m)+1$
%\State  $\tilde{\xxx}^{(j)}=\tilde{\xxx}_{i_{te}}$
%\If {$(i_{te}==h*m)$} 
%\State $h=h*2$ 
%\EndIf
%\EndIf
%\EndFor
%\end{algorithmic}
%\end{algorithm}

%\section{Section title of first appendix\label{app1}}

%Use \verb+\begin{verbatim}...\end{verbatim}+ for program codes without math. Use \verb+\begin{alltt}...\end{alltt}+ 

%\begin{lstlisting}[caption={Descriptive Caption Text},label=DescriptiveLabel]
%for i:=maxint to 0 do
%begin
%{ do nothing }
%end;
%Write('Case insensitive ');
%WritE('Pascal keywords.');
%\end{lstlisting}

%\nocite{*}% Show all bib entries - both cited and uncited; comment this line to view only cited bib entries;
%\bibliography{wileyNJD-VANCOUVER}%

%\section*{Author Biography}

%\begin{biography}{\includegraphics[width=66pt,height=86pt,draft]{empty}}{\textbf{Author Name.} This is sample author %\end{biography}

\end{document}